# On the Conditioned Exit Measures of Super Brownian Motion

By
Thomas S. Salisbury and John Verzani
*York University and the Fields Institute*
*CUNY – College of Staten Island*

*Summary:*

In this paper we present a martingale related to the exit measures of super Brownian motion. By changing measure with this martingale in the canonical way we have a new process associated with the conditioned exit measure. This measure is shown to be identical to a measure generated by a non-homogeneous branching particle system with immigration of mass. An application is given to the problem of conditioning the exit measure to hit a number of specified points on the boundary of a domain. The results are similar in flavor to the "immortal particle" picture of conditioned super Brownian motion but more general, as the change of measure is given by a martingale which need not arise from a single harmonic function.

## 1. Introduction

One may think about super Brownian motion heuristically as a measure-valued process induced by a branching particle system in which each particle diffuses as a Brownian motion and critical branching is performed. With this image in mind, Dynkin's exit measure for super Brownian motion from a domain $D$ is a measure on the boundary of $D$, supported on the set of points where the particles first exit the domain. Le Gall proved in [19] that the exit measure will "hit" a fixed point on the boundary of $D$, with positive probability, only in dimensions 1 and 2. In higher dimensions, a fixed point on the boundary will be hit by an exiting particle with probability 0. One may ask, if we condition the process to hit that point, what do the paths of the particles that hit the point look like? A naive guess is to say they are Brownian particles which are conditioned to hit the point, and so should look like $h$-transforms of Brownian motion. In dimension 2, where no conditioning is necessary, this is essentially correct. Although the first hit, as chronicled by the Brownian snake, is not an $h$-process, the typical path is. This has been exploited by Abraham [1].

A rough heuristic shows that one's intuition should be true in dimension $d \geq 3$. Let $z$ be a fixed point on $\partial D$, and for each $\epsilon > 0$ define $\Delta_\epsilon = \Delta_\epsilon(z) = \partial D \cap B(z, \epsilon)$ to be a ball on the boundary of $D$ centered at $z$ with radius $\epsilon$. Consider the "probability" $u_\epsilon(x)$ that the exit measure will charge the ball $\Delta_\epsilon$, under the excursion measure for the Brownian snake started from $x$. It is known (cf. Abraham and LeGall [2]) that if $D$ is regular enough, then $u_\epsilon$ will solve the equation $\Delta u = 4u^2$, with boundary condition that is infinite on $\Delta_\epsilon$ and 0 on its complement. Fix some point $x_0 \in D$. Then the ratio of functions $v_\epsilon(x) = u_\epsilon(x)/u_\epsilon(x_0)$ will solve the equation $\Delta v_\epsilon = 4v_\epsilon^2 u_\epsilon(x_0)$. As $\epsilon \to 0$, we will have that $u_\epsilon \to 0$, and thus $v_\epsilon$ should converge to a harmonic function. The boundary conditions suggest that the limit function is the harmonic function associated with Brownian motion conditioned to exit $D$ at $z$ – that is the Martin kernel $K_{x_0}(\cdot, z)$.

When this heuristic is formally inserted into the known formula for the Laplace transform, one finds that the exit measure conditioned to hit the point $z$ will look like an $h$-path transform of the exit measure. Let $D$ be a bounded domain in $\mathbb{R}^d$, and let $D_k \subseteq D$ be an increasing sequence of smooth subdomains with $D = \cup_k D_k$. Let $X^{D_k}$ denote the exit measure from $D_k$ under an excursion measure $\mathbb{N}_x$, let $\langle \cdot, \cdot \rangle$ denote integration and let $\phi$ be a test function defined on $\partial D_k$. Then formally,

$$\lim_{\epsilon \to 0} \mathbb{N}_{x_0}(\exp(-\langle X^{D_k}, \phi \rangle) \mid X^D(\Delta_\epsilon) > 0) = \mathbb{N}_{x_0}((\exp(-\langle X^{D_k}, \phi \rangle)\langle X^{D_k}, K_{x_0}(\cdot, z)\rangle).$$

In this case the "harmonic" function is now $H(\mu) = \langle \mu, K_{x_0}(\cdot, z) \rangle$. As is usual with $h$-path transforms, the change of measure is non-singular only for events prior to the exit time. Hence the appearance of $X^{D_k}$ in the statement.

In addition, there is a concrete representation of this $h$-transformed process. It is given in terms of an "immortal particle", or "backbone", consisting of a Brownian particle conditioned to exit $D$ at $z$. The conditioned super Brownian motion is obtained by allowing this particle to throw off mass in a uniform Poisson manner. This mass then evolves as with an (unconditioned) super Brownian motion.

Such $h$-transforms have been studied by Roelly-Coppoletta and Rouault [24], Evans and Perkins [16], and Overbeck [22], while the immortal particle representation originates with Evans [14]. For example, in [16] and in [24], super Brownian motion conditioned on non-extinction is shown to be described in terms of an $h$-path transform, where the harmonic function is now $H(\mu) = \langle \mu, 1 \rangle$. That is, if $X_t$ denotes the super Brownian motion in $\mathbb{R}^d$ at time $t$, then conditioning on non-extinction produces a Laplace functional of the form

$$\bar{\mathbb{N}}_x(\exp(-\langle X_t, \phi \rangle)) = \mathbb{N}_x((\exp(-\langle X_t, \phi \rangle))\langle X_t, 1\rangle).$$

The immortal particle representation in this case was obtained in [14]. For $h$ a harmonic function, [22] defines the transform by $H(\mu) = \langle \mu, h \rangle$. The immortal particle representation in that case was obtained in Overbeck [23].

Our interest is in more general conditionings, in which one conditions, not just upon one, but upon several specified boundary points being hit. In particular, we will show that if one fixes $n$ distinct points on the boundary of a bounded smooth domain $D$, and one conditions the exit measure to hit all of them, then the conditioned exit measure can be described in terms of an explicit martingale change of measure, or $h$-transform, the transforms of the preceding paragraph being a special case. We will then give a backbone representation of this $h$-transform, where now the backbone is a tree with $n$ leaves, each terminating at one of the $n$ boundary points. Our main purpose is to explicitly describe the evolution of this tree.

We note that in Etheridge [13], a different conditioning gives rise to a particle representation, where the backbone is also a tree. In Evans and O'Connell [15], a supercritical superprocess is given a particle representation in terms of a tree backbone throwing off subcritically evolving mass. In Dembo and Zeitouni [9], a large deviation result with some similarities to our problem gives rise to yet another tree structure.

We also note that steps have been taken towards a version of the Martin boundary, in the context of conditionings involving a single target (see Dynkin and Kuznetsov [12]). We conjecture that a more general Martin boundary theory would incorporate a modification of our basic conditionings, given in Section 6.

To give a more precise formulation of the results, we first describe a general class of $h$-transforms, which will include the ones we require as a particular case. Section 3 of the paper gives an analysis of this general class of transforms. To describe this class, we introduce some more notation. Let $g \geq 0$ be a solution to the non-linear equation $\frac{1}{2}\Delta u = 2u^2$ in $D$. Let $L_{4g}$ be the generator of Brownian motion killed at rate $4g$, and denote the law of the corresponding process by $E^{4g}$. Suppose we have $n$ positive solutions in $D$ to the linear equation $L_{4g}v = 0$, labeled $v^1, \ldots, v^n$. Let $U^{4g}$ be the potential operator for the generator $L_{4g}$ in $D$ and recursively define a family of functions $v^A$ indexed by nonempty subsets $A$ of $N = \{1, \ldots n\}$, as follows

$$v^A = \begin{cases} v^i & A = \{i\}, \\ 2\sum_{\substack{B \subseteq A \\ \emptyset, A \neq B}} U^{4g}(v^B v^{A \setminus B}) & |A| \geq 2. \end{cases}$$

For $\emptyset \neq A \subseteq N$ define

$$M_k^A = \sum_{\sigma \in \mathcal{P}(A)} \exp -\langle X^k, g \rangle \prod_{C \in \sigma} \langle X^k, v^C \rangle,$$

where $\mathcal{P}(A)$ is the set of partitions of $A$.

It turns out that $\{M_k^A\}_k$ forms a family of martingales for the exit measure with $\sigma$-fields $\mathcal{F}_k$ generated by the superprocess paths before they exit $D_k$.

**Theorem** (Theorem 3.1). *For each nonempty subset of $A \subset N$, $M_k^A$ is an $\mathcal{F}_k$-martingale under $\mathbb{N}_x$ with $\mathbb{N}_x(M_k^A) = v^A(x)$.*

Changing measure by $M_k^A$ thus consistently defines the $h$-transform

$$\mathbb{M}_x(\exp -\langle X^k, \phi \rangle) = \frac{1}{v^N(x)} \mathbb{N}_x(\exp(-\langle X^k, \phi \rangle) M_k^A).$$

Such a transform appears when one conditions the exit measure to hit $n$ distinct points on the boundary of the domain, as in the following result.

**Theorem** (Theorem 5.6). *Let $D$ be a bounded Lipschitz domain in $\mathbb{R}^d$, $d \geq 4$, and let $\{z_i\}_{i=1}^n$ be $n$ distinct points on the boundary of $D$. Set $g = 0$ and let $v^i = K_x(\cdot, z_i)$ be the Martin kernel with base point, $x$ and pole at $z_i$. Then*

$$\lim_{\epsilon \to 0} \mathbb{N}_x \left( \exp(-\langle X^k, \phi \rangle) \mid \prod_{i=1}^n X^D(\Delta_\epsilon(z_i)) > 0 \right) = \frac{1}{v^N(x)} \mathbb{N}_x(\exp(-\langle X^k, \phi \rangle) M_k^N).$$

Section 5 is largely devoted to proving this theorem. The arguments are based on an analysis of the asymptotics of small solutions to certain non-linear PDEs. Many of the basic estimates draw on the work of Abraham and LeGall [2]. In Section 6, we give a related

**preprint:** August 28, 1998

conjecture concerning conditioning super Brownian motion to hit the boundary in precisely $n$ points.

Having given an analytic formula for the conditioned exit measures, in terms of the martingale change of measure by $M_k^N$, we turn in Section 4 to generalizing the immortal particle picture and representing $\mathbb{M}_x$ dynamically in terms of a backbone shedding mass. This representation works for the general class of processes described above, with the effect of $g$ being to prune off mass that would otherwise reach $\partial D$. To summarize, in Theorem 4.2 we construct a backbone, which consists of a branching process that eventually branches into $n$ particles. In the example above these individual particles move as suitable transforms of Brownian motion. The particles in some sense keep track of a number of different exit points. Eventually there are $n$ particles, each of which moves as a Brownian motion conditioned to exit $D$ at a certain point. On top of this backbone, mass is immigrated in a Poisson manner and evolves as an unconditioned (but possibly pruned) super Brownian motion, to give us exit measures on $\{\partial D_k\}$. In the case there is just one particle to hit, then the backbone is simply an $h$-transform of Brownian motion. The arguments are based on an extended Palm formula for super Brownian motion.

Conditioning the exit measure to hit a point $z$ only produces a transform of the above type in dimensions high enough that the probability of such an event is actually zero. In dimension 2, any boundary point is hit with positive probability, and one must seek for other types of backbone representations for the conditioned exit measures. We give such a representation (among others) in an accompanying paper, [26].

For clarity, we will carry out the proof of the backbone representation in a simplified setting, avoiding the use of historical superprocesses. In fact, the arguments given would yield the stronger version of the representation, at the cost of a more baroque notation. In section 7 we will sketch the modifications needed to derive the strengthened results.

Section 2 gathers together some basic facts which are used in the rest of the paper.

## 2. Preliminaries

**2.1. Notation.** For a set $A$, let $|A|$ denote its cardinality, and let $\mathcal{P}(A)$ denote the collection of partitions of $A$. Choose some arbitrary linear order $\prec$ on the set of finite subsets of the integers. For $A$ such a finite subset, and $\sigma \in \mathcal{P}(A)$, let $\sigma(j)$ be the $j$th element of $\sigma$ in this order. Thus for example,

$$\prod_{C \in \sigma} \langle X^k, v^C \rangle = \prod_{j=1}^{|\sigma|} \langle X^k, v^{\sigma(j)} \rangle.$$

We will switch between these notations according to which seems clearer.

**2.2. Set facts.** We begin with some lemmas about sets. We will use the convention that a sum over an empty set is 0.

**Lemma 2.1.** *Let $A \subseteq B \subseteq C$ be subsets of $\{1, 2, \ldots n\}$. Then*

$$\sum_{A \subseteq B \subseteq C} (-1)^{|B|} = (-1)^{|C|} \mathbb{1}_{A=C}.$$

*Proof.*

$$\sum_{A \subseteq B \subseteq C} (-1)^{|B|} = (-1)^{|A|} \sum_{B \subseteq C \setminus A} (-1)^{|B|}$$

$$= (-1)^{|A|} \sum_{i=0}^{|C \setminus A|} \binom{|C \setminus A|}{i} (-1)^i$$

$$= \begin{cases} (-1)^{|C|} & C = A \\ (-1)^{|A|}(1-1)^{|C \setminus A|} & C \setminus A \neq \emptyset. \end{cases}$$

□

The following is immediate.

**Lemma 2.2.** *Let $A$ be finite, and let $w_i \in \mathbb{R}$ for $i \in A$. Then*

$$\prod_{i \in A}(1 - w_i) = 1 + \sum_{\substack{C \subseteq A \\ \emptyset \neq C}} (-1)^{|C|} \Big(\prod_{i \in C} w_i\Big).$$

In this paper we use the letter $K$ to denote a generic non-trivial constant whose particular value may vary from line to line. If it is important, explicit dependencies on other values will be specified.

**2.3. Facts about conditioned diffusions.** First we recall some formulae for conditioned Brownian motion.

Let $B$ be $d$-dimensional Brownian motion started from $x$, under a probability measure $P_x$. Write $\tau_D$ for the first exit time of $B$ from $D$.

Let $g : D \to [0, \infty)$ be bounded on compact subsets of $D$, and set

$$L_g = \frac{1}{2}\Delta - g.$$

Let $\xi_t$ be a process which, under a probability law $P_x^g$, has the law of a diffusion with generator $L_g$ started at $x$ and killed upon leaving $D$. In other words, $\xi$ is a Brownian motion on $D$, killed at rate $g$. Write $\zeta$ for the lifetime of $\xi$. Then

$$E_x^g(\xi_t \in A, \zeta > t) = E_x(\exp -\int_0^t ds\, g(B_s), B_t \in A, \tau_D > t). \tag{2.1}$$

Let $U^g f(x) = \int_0^\infty P_x^g(f(\xi_t)\mathbf{1}_{\{\zeta > t\}})dt$ be the potential operator for $L_g$. If $g = 0$ we write $U$ for $U^g$. If $0 \leq u$ is $L_g$-superharmonic, then the law of the $u$-transform of $\xi$ is determined by the formula

$$P_x^{g,u}(\Phi(\xi)\mathbf{1}_{\{\zeta > t\}}) = \frac{1}{u(x)} P_x^g(\Phi(\xi)u(\xi_t)\mathbf{1}_{\{\zeta > t\}})$$

for $\Phi(\xi) \in \sigma\{\xi_s; s \leq t\}$. Assuming that $0 < u < \infty$ on $D$, this defines a diffusion on $D$. If $u$ is $L_g$-harmonic, then it dies only upon reaching $\partial D$. If for $f \geq 0$, $u = U^g f$ (that is, $u$ is a potential) then it dies in the interior of $D$, and in fact $P_x^{g,u}$ satisfies

$$P_x^{g,u}(\Phi(\xi)) = \frac{1}{u(x)} \int_0^\infty P_x^g(\Phi(\xi_{\leq t})f(\xi_t)\mathbf{1}_{\{\zeta > t\}})dt, \tag{2.2}$$

where $\xi_{\leq t}$ is the process $\xi$ killed at time $t$.

### 2.4. Facts about the Brownian snake.
Next we recall some useful facts about the Brownian snake. Refer to [7] or [10] for a general introduction to superprocesses.

The Brownian snake is a path-valued process, devised by Le Gall as a means to construct super Brownian motion without limiting procedures. Refer to [18] or [20] for the construction.

We use the standard notation $(W_s, \zeta_s)$ for the Brownian snake, and $\mathbb{N}_x$ for the excursion measure of the Brownian snake starting from the trivial path $(w, \zeta), \zeta = 0, w(0) = x$. Note that $W_s(\cdot)$ is constant on $[\zeta_s, \infty)$, and $\zeta.$ has the distribution of a Brownian excursion under $\mathbb{N}_x$.

Super Brownian motion $X_t$ is defined as

$$\langle X_t, \phi \rangle = \int \phi(W_s(t)) \, dL_t(s),$$

where $L_t$ is the local time of $\zeta.$ at level $t$. Dynkin [11] introduced the exit measure $X^D$ associated with $X_t$. We follow Le Gall's snake-based definition of $X^D$ (see [20]) as

$$\langle X^D, \phi \rangle = \int \phi(W_s(\zeta_s)) \, dL^D(s),$$

where $L^D(\cdot)$ is an appropriate local time for $W_s(\zeta_s)$ on $\partial D$. (In section 7 we define these in a historical sense.)

We denote the range of the Brownian snake by $\mathcal{R}(W) = \{W_s(t) : 0 \le s \le \sigma, 0 \le t \le \zeta_s\}$ and the range inside $D$ by $\mathcal{R}^D(W) = \{W_s(t) : 0 \le s \le \sigma, 0 \le t \le \tau_D(W_s) \wedge \zeta_s\}$. There is an obvious inclusion between the range inside $D$ and the exit measures, given by

$$\{\langle X^D, \mathbf{1}_A \rangle > 0\} \subseteq \{\mathcal{R}^D(W) \cap A \ne \emptyset\}.$$

We make use of the following facts about the Brownian snake, which are contained in Le Gall [20].

A useful first-moment calculation is: (cf. [20] Proposition 3.3)

$$\mathbb{N}_x(\langle X^D, \phi \rangle) = E_x(\phi(B_{\tau_D})).$$

The following is an immediate consequence of Theorem 4.2 and its corollary in [20].

**Lemma 2.3.** *Let $g$ be a solution to $\Delta g = 4g^2$ in $D$, and let $\{D_k\}$ be an increasing sequence of smooth subdomains of $D$. Then for each $k$,*

$$\mathbb{N}_x(1 - \exp -\langle X^{D_k}, g \rangle) = g(x).$$

*Remark* 2.4. We make use of a generalization of this, replacing Brownian motion by Brownian motion killed at a certain rate.

Let $\mathcal{F}_k = \mathcal{F}_{D_k}$ be the $\sigma$-field of events determined by the superprocess killed upon exiting $D_k$. See [11] for a formal definition, or refer to section 7, where we give a definition in terms of the historical superprocess.

Dynkin introduced a Markov property for the exit measures in [11]. In our context, the Markov property is established in [21]. The next result gives it in the form we will use it:

**Lemma 2.5.**
$$\mathbb{N}_x(\exp -\langle X^D, \phi \rangle \mid \mathcal{F}_k) = \exp -\langle X^{D_k}, \mathbb{N}.(1 - \exp -\langle X^D, \phi \rangle) \rangle.$$

We use the following notation, where $B_s$ denotes a path in $D$ whose definition will be clear from the context:

$$e_\phi^D = e^D(\phi) = \exp -\langle X^D, \phi \rangle,$$

$$\mathcal{N}_t(e_\phi^D) = \mathcal{N}_t(e_\phi^D, B) = \exp - \int_0^t ds\, 4\mathbb{N}_{B_s}(1 - e_\phi^D).$$

The Palm formula for the Brownian snake takes the form: (cf. [20], Proposition 4.1)

$$\mathbb{N}_x(\langle X^D, \phi \rangle e_\psi^D) = E_x(\phi(B_{\tau_D})\mathcal{N}_{\tau_D}(e_\psi^D)). \tag{2.3}$$

We will make use of the following extension to the basic Palm formula. See [8] for a general discussion of this type of Palm formula.

**Lemma 2.6.** *Let $N = \{1, 2, \ldots n\}$, $n \geq 2$. Let $D$ be a domain, and let $B$ be a Brownian motion in $D$ with exit time $\tau$. Let $\{\psi_i\}$ be a family of measurable functions. Then*

$$\mathbb{N}_x\Big(e_\phi \prod_{i \in N}\langle X^D, \psi_i\rangle\Big)$$

$$= \frac{1}{2} \sum_{\substack{M \subseteq N \\ \emptyset, N \neq M}} E_x\Big(4 \int_0^\tau dt\, \mathcal{N}_t(e_\phi) \mathbb{N}_{B_t}\Big(e_\phi \prod_{i \in M}\langle X^D, \psi_i\rangle\Big)\mathbb{N}_{B_t}\Big(e_\phi \prod_{i \in N\setminus M}\langle X^D, \psi_i\rangle\Big)\Big).$$

*Proof.* Observe first that, by monotone convergence, it suffices to prove the result in the case that $D$ and the $\psi_i$ are bounded, and $\phi$ is bounded away from 0. So assume this.

Let $N^* = \{2, 3 \ldots, n\}$. Let $D_\lambda$ denote the derivative in $\lambda$. Then

$$\mathbb{N}_x\Big(e_\phi \prod_{i \in N}\langle X^D, \psi_i\rangle\Big) = \mathbb{N}_x\Big(\langle X^D, \psi_1\rangle e_\phi \prod_{i \in N^*}\langle X^D, \psi_i\rangle\Big) \tag{2.4}$$

$$= (-1)^{n-1} D_{\lambda_2} \cdots D_{\lambda_n} \mathbb{N}_x\Big(\langle X^D, \psi_1\rangle e(\phi + \sum_{i=2}^n \lambda_i \psi_i)\Big)\Big|_{\lambda_2 = \cdots = \lambda_n = 0}$$

$$= (-1)^{n-1} D_{\lambda_2} \cdots D_{\lambda_n} E_x\Big(\psi_1(B_\tau)\mathcal{N}_\tau\Big(e(\phi + \sum_{i=2}^n \lambda_i \psi_i)\Big)\Big)\Big|_{\lambda_2 = \cdots = \lambda_n = 0} \tag{2.5}$$

$$= \sum_{\sigma \in \mathcal{P}(N^*)} E_x\Big(\psi_1(B_\tau)\mathcal{N}_\tau(e_\phi)\prod_{j=1}^{|\sigma|}\Big(4\int_0^\tau dt\, \mathbb{N}_{B_t}\Big(e_\phi \prod_{i \in \sigma(j)}\langle X^D, \psi_i\rangle\Big)\Big)\Big).$$

Line (2.5) follows from the basic Palm formula. This last line uses the notation $\mathcal{P}(A)$ for all the partitions of a set $A$, and follows from a simple induction argument for the derivatives. Because $D$ is bounded, differentiation under the integral sign is easy to justify using dominated convergence, once we establish that $\mathbb{N}_x(e_\phi \prod_N\langle X^D, \psi_i\rangle)$ is bounded. But since the $\psi_i$ are bounded, and $\phi$ is bounded from 0, in fact $e_\phi \prod_{N^*}\langle X^D, \psi_i\rangle$ is itself bounded, so that this follows from (2.3) and the boundedness of $\psi_1$.

We use the notations $a\hat{b}c\ldots$ to denote all the elements $abc\ldots$ but the $b$th one. Let $\mathcal{G}_t$ be the filtration of $B_t$. Then by the Markov property for Brownian motion, one has that

$$\mathbb{N}_x\Big(e_\phi \prod_{i \in N}\langle X^D, \psi_i\rangle\Big)$$

$$= \sum_{\sigma \in \mathcal{P}(N^*)} E_x\Big(\psi_1(B_\tau)\mathcal{N}_\tau(e_\phi) \int_0^\tau \cdots \int_0^\tau dt_1 \cdots dt_{|\sigma|} \prod_{j=1}^{|\sigma|} 4\mathbb{N}_{B(t_j)}(e_\phi \prod_{i \in \sigma(j)} \langle X^D, \psi_i \rangle)\Big)$$

$$= \sum_{\sigma \in \mathcal{P}(N^*)} E_x\Big(\psi_1(B_\tau)\mathcal{N}_\tau(e_\phi) \sum_{k=1}^{|\sigma|} \int_0^\tau dt_k\, 4\mathbb{N}_{B(t_k)}(e_\phi \prod_{i \in \sigma(k)} \langle X^D, \psi_i \rangle)$$
$$\times \int_{t_k}^\tau \cdots \int_{t_k}^\tau dt_1 \cdots \hat{dt}_k \cdots dt_{|\sigma|} \prod_{j \neq k} 4\mathbb{N}_{B(t_j)}(e_\phi \prod_{i \in \sigma(j)} \langle X^D, \psi_i \rangle)\Big)$$

$$= \sum_{\sigma \in \mathcal{P}(N^*)} \sum_{k=1}^{|\sigma|} E_x\Big(4 \int_0^\tau dt_k\, \mathbb{N}_{B(t_k)}(e_\phi \prod_{i \in \sigma(k)} \langle X^D, \psi_i \rangle) E_x\Big(\psi_1(B_\tau)\mathcal{N}_\tau(e_\phi)$$
$$\times \int_{t_k}^\tau \cdots \int_{t_k}^\tau dt_1 \cdots \hat{dt}_k \cdots dt_{|\sigma|} \prod_{j \neq k} 4\mathbb{N}_{B(t_j)}(e_\phi \prod_{i \in \sigma(j)} \langle X^D, \psi_i \rangle) \mid \mathcal{G}_{t_k}\Big)\Big)$$

$$= \sum_{\emptyset \neq M \subseteq N^*} E_x\Big(4 \int_0^\tau dt\, \mathbb{N}_{B_t}(e_\phi \prod_{i \in M} \langle X^D, \psi_i \rangle) \mathcal{N}_t(e_\phi)$$
$$\times \sum_{\gamma \in \mathcal{P}(N^* \setminus M)} E_{B_t}\Big(\psi_1(B_\tau)\mathcal{N}_\tau(e_\phi) \prod_{j=1}^{|\gamma|} (4 \int_0^\tau dt\, \mathbb{N}_{B_t}(e_\phi \prod_{i \in \gamma(j)} \langle X^D, \psi_i \rangle))\Big)\Big)$$

$$= \sum_{\emptyset \neq M \subseteq N^*} E_x\Big(4 \int_0^\tau dt\, \mathcal{N}_t(e_\phi) \mathbb{N}_{B_t}(e_\phi \prod_{i \in M} \langle X^D, \psi_i \rangle) \mathbb{N}_{B_t}(e_\phi \prod_{i \in N \setminus M} \langle X^D, \psi_i \rangle)\Big)$$

$$= \frac{1}{2} \sum_{\substack{M \subseteq N \\ \emptyset, N \neq M}} E_x\Big(4 \int_0^\tau dt\, \mathcal{N}_t(e_\phi) \mathbb{N}_{B_t}(e_\phi \prod_{i \in M} \langle X^D, \psi_i \rangle) \mathbb{N}_{B_t}(e_\phi \prod_{i \in N \setminus M} \langle X^D, \psi_i \rangle)\Big).$$

Note that in order to rewrite the set of partitions of $N^*$, we designate one of the elements of a given partition as $M$, so that the rest forms a partition of $N \setminus M$. This isn't one-to-one, which is taken into account in the counting. The last factor of $1/2$ occurs because the summation in the last line counts everything twice. □

We use the extended Palm formula to show an exponential bound on the moments of the exit measure.

**Lemma 2.7.** *Let $D$ be a domain in $\mathbb{R}^d$ satisfying $\sup_{x \in D} E_x(\tau_D) < \infty$, where $\tau_D$ is the exit time from $D$ for Brownian motion. Then there exists $\lambda > 0$ such that*
$$\sup_{x \in D} \mathbb{N}_x(\exp \lambda \langle X^D, 1 \rangle - 1) < \infty.$$

*Remark* 2.8. A bounded domain $D$ in $\mathbb{R}^d$ will satisfy $\sup_D E_x(\tau_D) < \infty$.

*Proof.* It is enough to show that $\mathbb{N}_x(\langle X^D, 1 \rangle^n) \leq n! M^n$ for some $M < \infty$. Our proof follows part of the proof of Lemma 3.1 in Serlet [27]. Set $c_n = \sup_D \mathbb{N}_x(\langle X^D, 1 \rangle^n)$. Then $c_1 = \sup_D E_x(\tau_D) < \infty$ by assumption. By Lemma 2.6, with $\phi = 0$ we have the immediate recursion relation
$$c_1 = c_1$$

$$c_n \leq 2c_1 \sum_{j=1}^{n-1} \binom{n}{j} c_{n-j} c_j, \quad n \geq 2.$$

Let $a_{n+1} = c_1(2c_1^2)^n (2n)!/n!$, $a_1 = c_1$. Then one can easily verify, by induction, that the sequence $a_n$ satisfies the combinatorial identity

$$a_n = 2c_1 \sum_{j=1}^{n-1} \binom{n}{j} a_{n-j} a_j.$$

Since $c_1 = a_1 < \infty$, we get by induction that $c_n \leq a_n$. Stirling's formula applied to $a_n$ shows that some $M < \infty$ exists. □

## 3. A GENERAL CLASS OF MARTINGALES

We now present a general class of martingales related to the exit measure.

Let $D$ be a domain in $\mathbb{R}^d$, $d \geq 1$ and let $g \geq 0$ be a solution in $D$ to the non-linear equation

$$\frac{1}{2}\Delta u = 2u^2. \tag{3.1}$$

In particular, $g$ satisfies the conclusion of Lemma 2.3. Readers may feel free to take $g = 0$. In fact, the only application given in this paper requiring a non-zero $g$ is that of Section 6. Another reason we opt to work with a general $g$ is for consistency with [26], in which such $g$'s play an essential role.

Recall that $E^{4g}$ denotes the law of Brownian motion killed at rate $4g$, and that $L_{4g}$ is the generator of this process. That is

$$L_{4g} w = \frac{1}{2}\Delta w - 4gw.$$

Suppose we have $n$ positive solutions in $D$ to the linear equation $L_{4g} v = 0$, labeled $v^1, \ldots, v^n$. Recall that $U^{4g}$ is the potential operator for the generator $L_{4g}$ in $D$, and recursively define a family of functions $v^A$, for $\emptyset \neq A \subseteq N = \{1, \ldots n\}$, as follows:

$$v^A = \begin{cases} v^i & A = \{i\}, \\ 2\sum_{\substack{B \subseteq A \\ \emptyset, A \neq B}} U^{4g}(v^B v^{A\setminus B}) & |A| \geq 2. \end{cases} \tag{3.2}$$

Note that $v^A$ is either finite everywhere on $D$, or $v^A \equiv \infty$.

Let $D_k \Uparrow D$ be an increasing sequence of bounded, smooth subdomains and for each $k$ let $X^k$ be the exit measure from $D_k$, and $e_\phi^k = \exp -\langle X^k, \phi \rangle$. Let $\tau_k$ denote the exit time of a path from $D_k$. For $\emptyset \neq A \subseteq N$ define

$$M_k^A = \sum_{\sigma \in \mathcal{P}(A)} \exp(-\langle X^k, g \rangle) \prod_{C \in \sigma} \langle X^k, v^C \rangle. \tag{3.3}$$

It turns out that $\{M_k^A\}_k$ forms a family of martingales for the exit measure.

**Theorem 3.1.** *For each nonempty subset of $A \subset N$, $M_k^A$ is an $\mathcal{F}_k$-martingale under $\mathbb{N}_x$, provided $v^A < \infty$.*

*Remark* 3.2. Because $M_k^A$ is a martingale we can define a martingale change of measure in a canonical way as follows. For $\Phi_k$ a $\mathcal{F}_k$-measurable function, we set

$$\mathbb{M}_x(\Phi_k) = \frac{1}{v^A(x)} \mathbb{N}_x(\Phi_k M_k^A). \tag{3.4}$$

*Remark* 3.3. Letting $g = 0$ be the trivial solution to (3.1) and $v = 1$ (or $v = h$ with $h$ harmonic) we recover the well-known fact that $\langle X^k, 1 \rangle$, (or $\langle X^k, h \rangle$) is a martingale.

*Proof.* We first establish the following lemma

**Lemma 3.4.** *If $v^A < \infty$, then $\mathbb{N}_x(M_k^A) = v^A(x)$.*

*Proof.* The proof relies on the extended Palm formula. We induct on the size of $A$. First consider the case $|A| = 1$.

$$\begin{aligned}
\mathbb{N}_x(M_k^A) &= \mathbb{N}_x(\exp(-\langle X^k, g \rangle) \langle X^k, v^A \rangle) \\
&= E_x(v^A(\xi_{\tau_k}) \exp - \int_0^{\tau_k} ds\, 4 \mathbb{N}_{\xi_s}(1 - \exp - \langle X^k, g \rangle)) \quad &\text{(Palm formula)} \\
&= E_x(v^A(\xi_{\tau_k}) \exp - \int_0^{\tau_k} ds\, 4g(\xi_s)) \quad &\text{(Lemma (2.3))} \\
&= E_x^{4g}(v^A(\xi_{\tau_k}), \zeta > \tau_k) = v^A(x).
\end{aligned}$$

Next, suppose $|A| = n$ and that the lemma holds for all smaller sets. Then

$$v^A(x) = 2 \sum_{\substack{B \subseteq A \\ 0 < |B| < |A|}} U^{4g}(v^B v^{A \setminus B})(x)$$

$$= E_x^{4g}\left(v^A(\xi_{\tau_k}), \zeta < \tau_k\right) + 2 \sum_{\substack{B \subseteq A \\ 0 < |B| < |A|}} E_x\left(\int_0^{\tau_k} dt \exp(-\int_0^t ds\, 4g(\xi_s)) v^B(\xi_t) v^{A \setminus B}(\xi_t)\right).$$

Denote the first of the above terms by $I$, and the remainder by $II$. As in the case $|A| = 1$, $I = \mathbb{N}_x(e_g^k \langle X^k, v^A \rangle)$. Further

$$II = 2 \sum_{\substack{B \subseteq A \\ 0 < |B| < |A|}} E_x(\int_0^{\tau_k} dt\, \mathcal{N}_t(e_g^k) \mathbb{N}_{\xi_t}(M_k^B) \mathbb{N}_{\xi_t}(M_k^{A \setminus B})) \quad \text{(induction)}$$

$$= 2 \sum_{\substack{B \subseteq A \\ 0 < |B| < |A|}} E_x\left(\int_0^{\tau_k} dt\, \mathcal{N}_t(e_g^k) \sum_{\sigma \in \mathcal{P}(B)} \mathbb{N}_{\xi_t}(e_g^k \prod_{i=1}^{|\sigma|} \langle X^k, v^{\sigma(i)} \rangle) \right.$$
$$\left. \times \sum_{\gamma \in \mathcal{P}(A \setminus B)} \mathbb{N}_{\xi_t}(e_g^k \prod_{j=1}^{|\gamma|} \langle X^k, v^{\gamma(j)} \rangle)\right) \quad \text{(Definition of } M_k^A\text{)}$$

$$= \sum_{\substack{B \subseteq A \\ 0 < |B| < |A|}} \sum_{\sigma \in \mathcal{P}(B)} \sum_{\gamma \in \mathcal{P}(A \setminus B)} 2 E_x\left(\int_0^{\tau_k} dt\, \mathcal{N}_t(e_g^k) \right.$$
$$\left. \times \mathbb{N}_{\xi_t}(e_g^k \prod_{i=1}^{|\sigma|} \langle X^k, v^{\sigma(i)} \rangle) \mathbb{N}_{\xi_t}(e_g^k \prod_{j=1}^{|\gamma|} \langle X^k, v^{\gamma(j)} \rangle)\right)$$

$$= \sum_{\sigma \in \mathcal{P}(A)} \sum_{\substack{C \subseteq \{1,\ldots,|\sigma|\} \\ 0 < |C| < |\sigma|}} 2 E_x \Big( \int_0^{\tau_k} dt \, \mathcal{N}_t(e_g^k)$$

$$\times \mathbb{N}_{\xi_t}(e_g^k \prod_{i \in C} \langle X^k, v^{\sigma(i)} \rangle) \mathbb{N}_{\xi_t}(e_g^k \prod_{j \in C^c} \langle X^k, v^{\sigma(j)} \rangle) \Big)$$

$$= \sum_{\substack{\sigma \in \mathcal{P}(A) \\ |\sigma| > 1}} \mathbb{N}_x(e_g^k \prod_{i=1}^{|\sigma|} \langle X^k, v^{\sigma(i)} \rangle). \hspace{2cm} \text{(Palm formula)}$$

Thus

$$v^A(x) = I + II = \sum_{\sigma \in \mathcal{P}(A)} \mathbb{N}_x(e_g^k \prod_{i=1}^{|\sigma|} \langle X^k, v^{\sigma(i)} \rangle) = \mathbb{N}_x(M_k^A).$$

$\square$

Next we establish another lemma which helps us identify $\mathbb{N}_x(M_{k+1}^A \mid \mathcal{F}_k)$.

**Lemma 3.5.**

$$\mathbb{N}_x(M_{k+1}^A \mid \mathcal{F}_k) = e_g^k \sum_{\sigma \in \mathcal{P}(A)} \sum_{\gamma \in \mathcal{P}(\{1,\ldots,|\sigma|\})} \prod_{i=1}^{|\gamma|} \langle X^k, \mathbb{N}.(e_g^{k+1} \prod_{j \in \gamma(i)} \langle X^{k+1}, v^{\sigma(j)} \rangle) \rangle$$

*Proof.* This lemma follows from the strong Markov property satisfied by the exit measures (Lemma 2.5) and differentiation. Differentiation under the integral sign can be justified as before.

$\mathbb{N}_x(M_{k+1}^A \mid \mathcal{F}_k)$

$$= \sum_{\sigma \in \mathcal{P}(A)} (-1)^{|\sigma|} D_{\lambda_1} \cdots D_{\lambda_{|\sigma|}} \mathbb{N}_x \Big( \exp -\langle X^{k+1}, g + \sum_{i=1}^{|\sigma|} \lambda_i v^{\sigma(i)} \rangle \mid \mathcal{F}_k \Big) \Big|_{\lambda_1 = \cdots = \lambda_{|\sigma|} = 0}$$

$$= \sum_{\sigma \in \mathcal{P}(A)} (-1)^{|\sigma|} D_{\lambda_1} \cdots D_{\lambda_{|\sigma|}} \exp -\langle X^k,$$

(3.5)

$$\mathbb{N}. \Big( 1 - \exp -\langle X^{k+1}, g + \sum_{i=1}^{|\sigma|} \lambda_i v^{\sigma(i)} \rangle \Big) \rangle \Big|_{\lambda_1 = \cdots = \lambda_{|\sigma|} = 0}$$

$$= \sum_{\sigma \in \mathcal{P}(A)} \exp -\langle X^k, \mathbb{N}.(1 - \exp -\langle X^{k+1}, g \rangle) \rangle$$

(3.6)

$$\times \sum_{\gamma \in \mathcal{P}(\{1,\ldots,|\sigma|\})} \prod_{i=1}^{|\gamma|} \langle X^k, \mathbb{N}.(\exp(-\langle X^{k+1}, g \rangle) \prod_{j \in \gamma(i)} \langle X^{k+1}, v^{\sigma(j)} \rangle) \rangle.$$

$$= e_g^k \sum_{\sigma \in \mathcal{P}(A)} \sum_{\gamma \in \mathcal{P}(\{1,\ldots,|\sigma|\})} \prod_{i=1}^{|\gamma|} \langle X^k, \mathbb{N}.(e_g^{k+1} \prod_{j \in \gamma(i)} \langle X^{k+1}, v^{\sigma(j)} \rangle) \rangle.$$

Line (3.5) follows from the Markov property of the exit measures, and (3.6) is from the chain rule of calculus. $\square$

Finally, we prove the theorem. We have, by rearranging the following sum, that

$$M_k^A = e_g^k \sum_{\sigma \in \mathcal{P}(A)} \prod_{i=1}^{|\sigma|} \langle X^k, v^{\sigma(i)} \rangle$$

$$= e_g^k \sum_{\sigma \in \mathcal{P}(A)} \prod_{i=1}^{|\sigma|} \langle X^k, \mathbb{N}.(M_{k+1}^{\sigma(i)}) \rangle \tag{3.7}$$

$$= e_g^k \sum_{\sigma \in \mathcal{P}(A)} \prod_{i=1}^{|\sigma|} \sum_{\gamma \in \mathcal{P}(\sigma(i))} \langle X^k, \mathbb{N}.(e_g^{k+1} \prod_{j=1}^{|\gamma|} \langle X^{k+1}, v^{\gamma(j)} \rangle) \rangle$$

$$= e_g^k \sum_{\gamma \in \mathcal{P}(A)} \sum_{\sigma \in \mathcal{P}(\{1,\cdots,|\gamma|\})} \prod_{i=1}^{|\sigma|} \langle X^k, \mathbb{N}.(e_g^{k+1} \prod_{j \in \sigma(i)} \langle X^{k+1}, v^{\gamma(j)} \rangle) \rangle \tag{3.8}$$

$$= \mathbb{N}(M_{k+1}^A \mid \mathcal{F}_k). \tag{3.9}$$

Line (3.7) follows from Lemma 3.4, whereas line (3.8) follows from rearranging the terms in the sums. Line (3.9) is from Lemma 3.5. □

## 4. A BRANCHING PARTICLE DESCRIPTION

We now show that changing measure via $M_k^N$ is equivalent to taking a branching particle backbone process with immigration of mass along the paths of the particles. To formulate this, we define two measures, $\mathbb{M}_x$ and $\bar{\mathbb{N}}_x$. Recall that the former is defined by

$$\frac{d\mathbb{M}_x}{d\mathbb{N}_x}\Big|_{\mathcal{F}_k} = \frac{1}{v^N(x)} M_k^N,$$

so that

$$\mathbb{M}_x(\exp -\langle X^k, \phi \rangle) = \frac{1}{v^N(x)} \mathbb{N}_x(\exp(-\langle X^k, \phi \rangle) M_k^N).$$

To define $\bar{\mathbb{N}}_x$, we first construct the branching process which will form the backbone of our new exit measures. Let $D$, $D_k$, $N$, $g$ and $\{v^A\}$ be as before. By (3.2) and the linearity of $U^{4g}$, we have that $v^A$ is an $L_{4g}$-potential for $|A| \geq 2$ (provided it is finite). For $|A| = 1$ it is $L_{4g}$-harmonic. The branching-particle process we desire will start with a single $v^N$-particle, that is, a particle moving as a $v^N$-transform of the process with generator $L_{4g}$. If $|N| = 1$, then because $v^N$ is $L_{4g}$-harmonic, the particle dies on the boundary of $D$. If $|N| \geq 2$ then $v^N$ is an $L_{4g}$-potential and the particle dies in the interior of $D$. At its death, it splits into two independent particles, a $v^A$-particle and a $v^{N \setminus A}$ particle, with position-dependent probability given by

$$p(A, N)(y) = \frac{(v^A v^{N \setminus A})(y)}{\sum_{\emptyset \neq B \neq N} (v^B v^{N \setminus B})(y)}. \tag{4.1}$$

The $v^A$ particle tracks the functions $\{v^i\}_{i \in A}$ while the $v^{N \setminus A}$ particle tracks the remaining functions. This pattern then repeats for each new particle. For example, the $v^A$ particle dies on the boundary of $D$ if $A$ is a singleton, and otherwise it dies in the interior of $D$, giving birth to two new, independent particles as before. Since the lifetime of each particle is finite almost surely, eventually all of the particles will have died out. In total there will be $2n - 1$ particles of which $n$ die on the boundary of $D$. The remaining $n - 1$ will die in the interior of $D$.

Let $n_t$ denote the number of particles alive at time $t$. Label them with $1 \leq i \leq n_t$, and for each one set $x^i(s), 0 \leq s \leq t$ to be the history (including the ancestors' history) of the individual particle up until time $t$. Define measure-valued branching processes as follows

$$\Upsilon_t(dx) = \sum_{i=1}^{n_t} \delta_{x^i(t)}(dx), \quad \Upsilon_t^k(dx) = \sum_{i=1}^{n_t} \mathbf{1}_{\tau_k(x^i) > t} \delta_{x^i(t)}(dx). \tag{4.2}$$

The process $\Upsilon_t^k$ puts a mass at each particle alive at time $t$ which hasn't already exited $D_k$. Without comment, these processes will be referred to in terms of the underlying particles although strictly speaking they are measures. Let $Q_x$ denote the law of $\Upsilon_t$.

We now specify how the mass $\Upsilon$ throws off evolves. The following calculation shows that $\exp-\langle X^k, g\rangle$ is an $\mathcal{F}_k$-martingale (although it has an infinite moment under $\mathbb{N}_x$), and follows by the Markov property and Lemma 2.3.

$$\mathbb{N}_x(\exp-\langle X^k, g\rangle \mid \mathcal{F}_{k-1}) = \exp-\langle X^{k-1}, \mathbb{N}.(1 - \exp-\langle X^k, g\rangle)\rangle$$
$$= \exp-\langle X^{k-1}, g\rangle.$$

Thus we can define a consistent measure $\tilde{\mathbb{N}}_x$ on the $\mathcal{F}_k$ measurable sets by

$$\tilde{\mathbb{N}}_x(\Phi_k) = \mathbb{N}_x(\Phi_k \exp-\langle X^k, g\rangle). \tag{4.3}$$

The extra factor prunes off mass that would otherwise get to $\partial D_k$. Though we will not need it in the rest of the argument, the following is a justification for this interpretation of $\tilde{\mathbb{N}}_x$. It is essentially a special case of Dawson's Girsanov formula (see section 10.1.2 of [7]).

**Lemma 4.1.** *Let $\tilde{\tilde{\mathbb{N}}}_x$ be the excursion law for the superprocess in $D$ based on the generator $L_{4g}$. Then $\tilde{\mathbb{N}}_x(1 - \exp-\langle X^k, \phi\rangle) = \tilde{\tilde{\mathbb{N}}}_x(1 - \exp-\langle X^k, \phi\rangle)$ for every $\phi \geq 0$.*

*Proof.* Let $\psi$ be the solution to

$$L_{4g} u = 2u^2 \text{ in } D_k$$
$$u = \phi \text{ on } \partial D_k.$$

It is then easily checked that $u = \psi + g$ is the solution to

$$\frac{1}{2}\Delta u = 2u^2 \text{ in } D_k$$
$$u = \phi + g \text{ on } \partial D_k,$$

Thus by Lemma 2.3,

$$\tilde{\mathbb{N}}_x(1 - \exp-\langle X^k, \phi\rangle) = \mathbb{N}_x\Big((1 - \exp-\langle X^k, \phi\rangle)\exp-\langle X^k, g\rangle\Big)$$
$$= \mathbb{N}_x(1 - \exp-\langle X^k, \phi + g\rangle) - \mathbb{N}_x(1 - \exp-\langle X^k, g\rangle)$$
$$= \Big(\psi(x) + g(x)\Big) - g(x) = \psi(x)$$
$$= \tilde{\tilde{\mathbb{N}}}_x(1 - \exp-\langle X^k, \phi\rangle).$$

This suffices. □

One thinks of mass being created continuously along the backbone, but only at countably many times will it actually survive, even instantaneously. At each such time, the mass created evolves like a superprocess with "law" some $\tilde{\mathbb{N}}_y$, and, if it survives long enough, produces a contribution to the exit measure. More properly, given the backbone $\Upsilon^k$, we form a Poisson random measure $N^k(d\mu)$ with intensity $\int_0^\infty dt \int 4\Upsilon_t^k(dy)\tilde{\mathbb{N}}_y(X^k \in d\mu)$. We then realize the exit measure under $\bar{\mathbb{N}}_x$ as $Y^k = \int \mu N^k(d\mu)$. A standard calculation shows that

$$\bar{\mathbb{N}}_x(\exp -\langle Y^k, \phi \rangle) = Q_x(\exp - \int_0^\infty dt\, 4\langle \Upsilon_t^k, \tilde{\mathbb{N}}.(1 - \exp -\langle X^k, \phi \rangle)\rangle). \tag{4.4}$$

**Theorem 4.2.** *If $v^N < \infty$, then $\mathbb{M}_x(\exp -\langle X^k, \phi \rangle) = \bar{\mathbb{N}}_x(\exp -\langle Y^k, \phi \rangle)$ for each $\phi \geq 0$.*

*Remark* 4.3. In section 7 we will sketch a proof that in fact, $\mathbb{M}_x = \bar{\mathbb{N}}_x$ on $\mathcal{F}_k$.

*Proof.* Recall that

$$\mathbb{M}_x(\exp -\langle X^k, \phi \rangle) = \frac{1}{v^N(x)}\mathbb{N}_x(e_\phi^k M_k^N).$$

For $A \subseteq N$, $\sigma \in \mathcal{P}(A)$ define $m_k^\sigma = e_g^k \prod_{B \in \sigma} \langle X^k, v^B \rangle$, so that $M_k^A = \sum_{\sigma \in \mathcal{P}(A)} m_k^\sigma$. For $B \subseteq A$ define

$$\sigma|_B = \begin{cases} \{B_1, \ldots, B_m\} & B = \cup_{j=1}^m B_j, \text{ each } B_j \in \sigma \\ \emptyset & \text{if no such decomposition exists.} \end{cases}$$

That is, if $\sigma$ can be restricted to be a partition of $B$, then define it as such, otherwise set $\sigma|_B$ to be the empty set.

We use $\xi$ for the canonical process with lifetime $\zeta$.

We first mention a consequence of the Palm formula. Let $\sigma \in \mathcal{P}(N)$. Then

$$\mathbb{N}_x(e_\phi^k m_k^\sigma) = \frac{1}{2} \sum_{\substack{C \subset \{1,\ldots,|\sigma|\} \\ 1 < |C| < |\sigma|}} E_x\Big(4 \int_0^{\tau_k} dt\, \mathcal{N}_t(e_{\phi+g}^k)\mathbb{N}_{\xi_t}(e_{\phi+g}^k \prod_{i \in C} \langle X^k, v^{\sigma(i)} \rangle)$$

$$\times \mathbb{N}_{\xi_t}(e_{\phi+g}^k \prod_{i \in \{1,\ldots,|\sigma|\}\setminus C} \langle X^k, v^{\sigma(i)} \rangle)\Big)$$

$$= \frac{1}{2} \sum_{\substack{A \subseteq N \\ \emptyset, N \neq A \\ \sigma|_A \neq \emptyset}} E_x(4 \int_0^{\tau_k} dt\, \mathcal{N}_t(e_{\phi+g}^k)\mathbb{N}_{\xi_t}(e_\phi^k m_k^{\sigma|_A})\mathbb{N}_{\xi_t}(e_\phi^k m_k^{\sigma|_{N\setminus A}})). \tag{4.5}$$

The branching processes $\Upsilon$ naturally partitions the $|N|$ points in the following manner. We eventually have $n$ particles that exit the domain $D$. Group together those that have a common ancestor that exits $D_k$. We will write $\Upsilon^k \sim \sigma$ for the relationship that the resulting partition is $\sigma$.

**Lemma 4.4.** *On the event $\{\Upsilon^k \sim \sigma\}$ one has*

$$Q_x(\exp - \int_0^\infty dt\, 4\langle \Upsilon_t^k, \tilde{\mathbb{N}}.(1 - \exp -\langle X^k, \phi \rangle)\rangle); \Upsilon^k \sim \sigma) = \frac{1}{v^N(x)}\mathbb{N}_x(e_\phi^k m_k^\sigma). \tag{4.6}$$

- *page 14* -

*Proof.* This will be established by induction on the size of $|\sigma|$. First we look at the case $|\sigma| = 1$, in other words, where there is no branching and $\sigma = \{N\}$. Since $\Upsilon_t^k$ is formed by a $v^N$-particle which exits $D_k$ before dying we have

$$Q_x(\exp - \int_0^\infty dt\, 4\langle \Upsilon_t^k, \tilde{\mathbb{N}}.(1 - \exp -\langle X^k, \phi\rangle)\rangle, \Upsilon^k \sim \sigma)$$

$$= E_x^{4g, v^N}(\exp - \int_0^{\tau_k} dt\, 4\tilde{\mathbb{N}}_{\xi_t}(1 - \exp -\langle X^k, \phi\rangle), \zeta > \tau_k)$$

$$= \frac{1}{v^N(x)} E_x^{4g}\left(v^N(\xi_{\tau_k}) \exp\left(-\int_0^{\tau_k} dt\, 4\mathbb{N}_{\xi_t}\left(\exp -\langle X^k, g\rangle(1 - \exp -\langle X^k, \phi\rangle)\right)\right), \zeta > \tau_k\right)$$

$$= \frac{1}{v^N(x)} E_x\left(v^N(\xi_{\tau_k}) \exp(-\int_0^{\tau_k} dt\, 4g(\xi_t)) \exp(\int_0^{\tau_k} dt\, 4\mathbb{N}_{\xi_t}(1 - \exp -\langle X^k, g\rangle))\right.$$
$$\left. \times \exp(-\int_0^{\tau_k} dt\, 4\mathbb{N}_{\xi_t}(1 - \exp -\langle X^k, \phi + g\rangle))\right)$$

$$= \frac{1}{v^N(x)} E_x\left(v^N(\xi_{\tau_k}) \exp(-\int_0^{\tau_k} dt\, 4g(\xi_t)) \exp(\int_0^{\tau_k} dt\, 4g(\xi_t))\right.$$
$$\left. \times \exp(-\int_0^{\tau_k} dt\, 4\mathbb{N}_{\xi_t}(1 - \exp -\langle X^k, \phi + g\rangle))\right)$$

$$= \frac{1}{v^N(x)} E_x(v^N(\xi_{\tau_k}) \mathcal{N}_{\tau_k}(e_{\phi+g}^k))$$

$$= \frac{1}{v^N(x)} \mathbb{N}_x(e_\phi^k e_g^k \langle X^k, v^N\rangle) \tag{4.7}$$

$$= \frac{1}{v^N(x)} \mathbb{N}_x(e_\phi^k m_k^\sigma).$$

Line (4.7) follows from the Palm formula with $\sigma = \{N\}$.

To establish (4.6) for all $\sigma \in \mathcal{P}(N)$, we assume the induction hypothesis is true for all $A \subseteq N$ and all partitions of $A$ of size $j$ or smaller, $j \geq 1$. Now consider $\sigma \in \mathcal{P}(N)$ of size $j + 1$.

Since $|\sigma| \geq 2$, we must have that $\Upsilon$ branches before it exits $D_k$. The first branch time of $\Upsilon$ is the lifetime $\zeta$ of the $v^N$-particle.

$$Q_x\left(\exp\left(-\int_0^\infty dt\, 4\langle \Upsilon_t^k, \tilde{\mathbb{N}}.(1 - \exp -\langle X^k, \phi\rangle)\rangle\right); \Upsilon^k \sim \sigma\right)$$

$$= Q_x\left(\exp\left(-\left(\int_0^\zeta + \int_\zeta^\infty\right) dt\, 4\langle \Upsilon_t^k, \tilde{\mathbb{N}}.(1 - \exp -\langle X^k, \phi\rangle)\rangle\right); \Upsilon^k \sim \sigma\right)$$

$$= E_x^{4g, v^N}\left(\mathbb{1}_{\zeta < \tau_k} \exp\left(-\int_0^\zeta dt\, 4\tilde{\mathbb{N}}_{\xi_t}(1 - \exp -\langle X^k, \phi\rangle)\right) \sum_{\substack{A \subseteq N \\ \emptyset, A \neq N \\ \sigma|_A \neq \emptyset}} p(A, N)(\xi_\zeta)\right.$$
$$\times Q_{\xi_\zeta}\left(\exp\left(-\int_0^\infty dt\, 4\langle \Upsilon_t^k, \tilde{\mathbb{N}}.(1 - \exp -\langle X^k, \phi\rangle)\rangle\right); \Upsilon^k \sim \sigma|_A\right)$$
$$\left. \times Q_{\xi_\zeta}\left(\exp\left(-\int_0^\infty dt\, 4\langle \Upsilon_t^k, \tilde{\mathbb{N}}.(1 - \exp -\langle X^k, \phi\rangle)\rangle\right); \Upsilon^k \sim \sigma|_{N \setminus A}\right)\right) \tag{4.8}$$

$$= E_x^{4g,v^N}\Bigl(\mathbf{1}_{\zeta<\tau_k} \exp\Bigl(-\int_0^\zeta dt\, 4\tilde{\mathbb{N}}_{\xi_t}(1-\exp-\langle X^k,\phi\rangle)\Bigr)\Bigr) \sum_{\substack{A\subset N \\ \emptyset,A\neq N \\ \sigma|_A\neq\emptyset}} p(A,N)(\xi_\zeta) \tag{4.9}$$

$$\times \frac{1}{v^A(\xi_\zeta)}\mathbb{N}_{\xi_\zeta}(e_\phi^k m_k^{\sigma|_A})\frac{1}{v^{N\setminus A}(\xi_\zeta)}\mathbb{N}_{\xi_\zeta}(e_\phi^k m_k^{\sigma|_{N\setminus A}}))$$

$$= \sum_{\substack{A\subset N \\ \emptyset,A\neq N \\ \sigma|_A\neq\emptyset}} E_x^{4g,v^N}\Bigl(\mathbf{1}_{\zeta<\tau_k}\exp\Bigl(-\int_0^\zeta dt\, 4\tilde{\mathbb{N}}_{\xi_t}(1-\exp-\langle X^k,\phi\rangle)\Bigr)\Bigr)$$

$$\times \frac{1}{\sum_{\substack{B\subset N \\ \emptyset,N\neq B}}(v^B v^{N\setminus B})(\xi_\zeta)}\mathbb{N}_{\xi_\zeta}(e_\phi^k m_k^{\sigma|_A})\mathbb{N}_{\xi_\zeta}(e_\phi^k m_k^{\sigma|_{N\setminus A}}))$$

$$= \sum_{\substack{A\subset N \\ \emptyset,A\neq N \\ \sigma|_A\neq\emptyset}} \frac{1}{v^N(x)} E_x^{4g}\Bigl(\int_0^{\tau_k} ds\, \mathbf{1}_{\zeta>s}\exp\Bigl(-\int_0^s dt\, 4\tilde{\mathbb{N}}_{\xi_t}(1-\exp-\langle X^k,\phi\rangle)\Bigr)\Bigr)$$

$$\times \frac{1}{\sum_{\substack{B\subset N \\ \emptyset,N\neq B}}(v^B v^{N\setminus B})(\xi_s)}\mathbb{N}_{\xi_s}(e_\phi^k m_k^{\sigma|_A})\mathbb{N}_{\xi_s}(e_\phi^k m_k^{\sigma|_{N\setminus A}})(2\sum_{\substack{B\subset N \\ \emptyset,N\neq B}}(v^B v^{N\setminus B})(\xi_s)))$$
$$\tag{4.10}$$

$$= 2\sum_{\substack{A\subset N \\ \emptyset,A\neq N \\ \sigma|_A\neq\emptyset}} \frac{1}{v^N(x)} E_x\Bigl(\int_0^{\tau_k} ds\, \exp\Bigl(-\int_0^s dt\, 4g(\xi_t)\Bigr)\exp\Bigl(\int_0^s dt\, 4\mathbb{N}_{\xi_t}(1-\exp-\langle X^k,g\rangle)\Bigr)$$

$$\times \exp\Bigl(-\int_0^s dt\, 4\mathbb{N}_{\xi_t}(1-\exp-\langle X^k,\phi+g\rangle)\Bigr)\mathbb{N}_{\xi_s}(e_\phi^k m_k^{\sigma|_A})\mathbb{N}_{\xi_s}(e_\phi^k m_k^{\sigma|_{N\setminus A}}))$$
$$\tag{4.11}$$

$$= 2\sum_{\substack{A\subset N \\ \emptyset,A\neq N \\ \sigma|_A\neq\emptyset}} \frac{1}{v^N(x)} E_x\Bigl(\int_0^{\tau_k} ds\, \mathcal{N}_s(e_{\phi+g}^k)\mathbb{N}_{\xi_s}(e_\phi^k m_k^{\sigma|_A})\mathbb{N}_{\xi_s}(e_\phi^k m_k^{\sigma|_{N\setminus A}}))$$

$$= \mathbb{N}_x(e_\phi^k m_k^\sigma). \tag{4.12}$$

Line (4.8) follows from the strong Markov property for the branching process $\Upsilon$ applied at the first branch time, and the description of the branching process having independent offspring. The next line, (4.9) then follows from the induction hypothesis. The formula (2.2) for a conditioned $u$-processes gives line (4.10) and (4.11) follows since under $E^{4g}$ the process is Brownian motion killed at rate $4g$. Finally, (4.12) comes from (4.5). □

By summing (4.6) over $\sigma \in \mathcal{P}(N)$ we get

$$Q_x(\exp-\int_0^\infty dt\, 4\langle \Upsilon_t^k, \tilde{\mathbb{N}}.(1-\exp-\langle X^k,\phi\rangle)\rangle)$$

$$= \frac{1}{v^N(x)}\mathbb{N}_x(e_\phi^k \sum_{\sigma\in\mathcal{P}(N)} m_k^\sigma) = \mathbb{N}_x(e_\phi^k M_k^N) = \mathbb{M}_x(e_\phi^k).$$

Using the definition of $\bar{\mathbb{N}}$ from (4.4) we get that
$$\bar{\mathbb{N}}_x(\exp -\langle Y^k, \phi \rangle) = \mathbb{M}_x(\exp -\langle X^k, \phi \rangle).$$
□

## 5. Conditioning the support of the exit measure to hit $n$ points on $\partial D$

In this section we investigate the exit measure when it is conditioned to charge $n$ small balls on the boundary of $D$. In the limit as the radius of the balls tends to 0, the conditioned process converges to one given by a martingale change of measure as in section 3.

Let $N = \{1, 2, \ldots, n\}$ and let $\{z_i\}_{i \in N}$ be a finite set of distinct points on the boundary of a bounded, Lipschitz domain $D$. Using the notation $B(x, \epsilon)$ for a ball centered at $x$ with radius $\epsilon$, define the sets
$$B^i_\epsilon = B(z_i, \epsilon),$$
$$\Delta^i_\epsilon = B^i_\epsilon \cap \partial D.$$
For a set $A \subseteq N$ set $B^A_\epsilon = \cup_A B^i_\epsilon$ and $\Delta^A_\epsilon = \cup_A \Delta^i_\epsilon$. Furthermore, for $x \in D$ define the functions
$$u^A_\epsilon(x) = \mathbb{N}_x(\sum_{i \in A} \langle X^D, \mathbf{1}_{\Delta(z_i, \epsilon)} \rangle > 0) = \mathbb{N}_x(\text{Hit at least one of } \Delta(z_i, \epsilon) \text{ while exiting}),$$
$$v^A_\epsilon(x) = \mathbb{N}_x(\prod_{i \in A} \langle X^D, \mathbf{1}_{\Delta(z_i, \epsilon)} \rangle > 0) = \mathbb{N}_x(\text{Hit all the } \Delta(z_i, \epsilon) \text{ while exiting}).$$

It follows from Theorem 7.1 of [2] that
$$\frac{1}{2} \Delta u^A_\epsilon = 2(u^A_\epsilon)^2$$
$$\lim_{y \to z \in \partial D} u^A_\epsilon(y) = \begin{cases} 0 & z \in \partial D \setminus \overline{B^A_\epsilon} \\ \infty & z \in B^A_\epsilon \end{cases}$$

A simple relationship between the two functions is summarized in this lemma which is proved by an inclusion exclusion argument.

**Lemma 5.1.** *We have the following relationships:*
$$v^A_\epsilon = -\sum_{\emptyset \neq B \subseteq A} (-1)^{|B|} u^B_\epsilon$$
$$u^A_\epsilon = -\sum_{\emptyset \neq B \subseteq A} (-1)^{|B|} v^B_\epsilon.$$

We now show that the family of functions $v^A_\epsilon$ satisfies the following recursive relations.

**Theorem 5.2.** *For $D$ a bounded Lipschitz domain, and $A$ and $v^A_\epsilon$ as above, the functions $v^A_\epsilon$ satisfy*

(a) $\quad v^A_\epsilon \in \mathcal{C}^2(D)$.

(b) $\quad \dfrac{1}{2} \Delta v^A_\epsilon = -2 \sum_{\substack{\emptyset \neq B, C \subseteq A \\ B \cup C = A}} (-1)^{|B \cap C|} v^B_\epsilon v^C_\epsilon$.

$$= 2\left(\left(2u_\epsilon^A + (-1)^{|A|}v_\epsilon^A\right)v_\epsilon^A - \sum_{\substack{B\cup C=A \\ \emptyset, A\neq B,C}} (-1)^{|B\cap C|}v_\epsilon^B v_\epsilon^C\right).$$

$$(c) \quad \lim_{y\to z\in\partial D} v_\epsilon^A(y) = \begin{cases} 0 & |A|\geq 2 \\ \begin{cases} 0 & z\in\partial D\setminus\overline{\Delta_\epsilon^A} \\ \infty & z\in\Delta_\epsilon^A \end{cases} & |A|=1. \end{cases}$$

*Proof.* For convenience, drop the $\epsilon$ subscript. When $|A|=1$ this is Theorem 7.1 of [2]. Let $\Theta$ denote surface measure on $\partial D$, which exists by the Lipschitz condition. The condition on $D$ in [2] is that it satisfy

$$\liminf_{r\to 0} \frac{\Theta(\Delta_r^a \cup \Delta_\epsilon^A)}{\Theta(\Delta_r^a)} > 0, \text{ for all } a\in\Delta_\epsilon^A.$$

For this problem, this is trivially true as the limit is 1.

For $|A|\geq 2$, we have by Lemma 5.1 and the relation $(1/2)\Delta u^A = 2(u^A)^2$, that

$$\frac{1}{2}\Delta v^A = -2\sum_{\emptyset\neq B\subseteq A} (-1)^{|B|}(u^B)^2$$

$$= -2\sum_{\emptyset\neq B\subseteq A} (-1)^{|B|}\left(\sum_{\emptyset\neq C\subseteq B} (-1)^{|C|}v^C\right)^2 \tag{5.1}$$

$$= -2\sum_{\emptyset\neq B\subseteq A} (-1)^{|B|} \sum_{\emptyset\neq C,D\subseteq B} (-1)^{|C|+|D|}v^C v^D$$

$$= -2\sum_{\emptyset\neq C,D\subseteq A} (-1)^{|C|+|D|}v^C v^D \sum_{C\cup D\subseteq B\subseteq A} (-1)^{|B|}$$

$$= -2\sum_{\emptyset\neq C,D\subseteq A} (-1)^{|A|+|C\cap D|}v^C v^D (-1)^{|A|}\mathbf{1}_{C\cup D=A} \tag{5.2}$$

$$= -2\left(2\sum_{\emptyset\neq C\subseteq A} (-1)^{|C|}v^C v^A - (-1)^{|A|}(v^A)^2\right) - 2\sum_{\substack{C\cup D=A \\ \emptyset, A\neq C,D}} (-1)^{|C\cap D|}v^C v^D$$

$$= -2\left(2(-u^A)v^A - (-1)^{|A|}(v^A)^2\right) - 2\sum_{\substack{C\cup D=A \\ \emptyset, A\neq C,D}} (-1)^{|C\cap D|}v^C v^D.$$

Here line (5.1) follows by Lemma 5.1, and line (5.2) by Lemma 2.1.

Finally, we remark that $v^A(y)\leq v^i(y)$ for $i\in A$. Fix $z\in\partial D$. For $\epsilon$ sufficiently small, because of the Lipschitz assumption on $D$, $z$ is not contained in $B_\epsilon^i$ for some $i\in A$. Thus as $y\to z$ we, have $0\leq v_\epsilon^A(y)\leq v_\epsilon^i(y)\to 0$ by the $|A|=1$ case. $\square$

**Theorem 5.3.** *Let $D$ be a bounded Lipschitz domain in $\mathbb{R}^d$ for $d\geq 4$. Let $\{z_1,\ldots z_n\}$ be distinct points on the boundary of $D$, $A\subseteq N = \{1,2,\ldots n\}$, and $x,y\in D$. Let $U$ be the potential operator for Brownian motion killed upon hitting the boundary of $D$, $K_x^D$ the Martin kernel for $D$. Then, the following limit exists*

$$\lim_{\epsilon\to 0} \frac{v_\epsilon^A(y)}{\prod_A v_\epsilon^i(x)} = \Psi_x^A(y).$$

*Furthermore, the limit satisfies*

$$\Psi_x^A(y) = \begin{cases} K_x^D(y, z_i) & A = \{i\}, \\ 2\sum_{\substack{B \subseteq A \\ A, \emptyset \neq B}} U(\Psi_x^B \Psi_x^{A \setminus B})(y) & |A| \geq 2. \end{cases}$$

*Proof.* Fix $\delta > 0$ such that $B(z_i, \delta) \cap B(z_j, \delta) = \emptyset$ if $i \neq j \in A$. Set $D_\epsilon^A = D \setminus \overline{\cup_A B(z_i, \epsilon)}$, and $\tau_\epsilon = \tau_\epsilon^A = \tau_{D_\epsilon^A}$.

The proof of this theorem relies on the following lemma.

**Lemma 5.4.** *Assume the conditions of Theorem 5.3, and fix an $x \in D$. There exist constants $K < \infty$ and $\epsilon_0 > 0$ such that for all $\epsilon < \epsilon_0$, $A \subseteq N$ and $y \in D_{4|A|\epsilon}^A$*

$$\frac{v_\epsilon^A(y)}{\prod_A v_\epsilon^i(x)} \leq K \sum_A K_x^D(y, z_i).$$

*In particular, the expression on the left is bounded in $y$ on any compact subset of $D$, uniformly in $\epsilon$.*

*Proof (of lemma).* We prove the lemma by induction on the size of $A$. First the case $A = \{i\}$. Set $z = z_i$, and $v = v^i$. This part simply extends the argument of the proof of Theorem 3.1 of [2] to the Lipschitz case.

Since $\Delta v_\epsilon = 4(v_\epsilon)^2$ in $D$ and $v_\epsilon$ vanishes on $\partial D \setminus \Delta(z, \epsilon)$ the Feynman-Kac formula gives

$$v_\epsilon(y) = P_y\Big(v_\epsilon(B(\tau_{2\epsilon})) \exp(-\int_0^{\tau_{2\epsilon}} 2v_\epsilon(B_s)ds)\Big)$$
$$\leq P_y\Big(v_\epsilon(B(\tau_{2\epsilon})), B(\tau_{2\epsilon}) \in \partial B(z, 2\epsilon)\Big). \tag{5.3}$$

But, from the definition of $v_\epsilon$ one has the bound

$$v_\epsilon(w) = \mathbb{N}_w(X^D(B(z, \epsilon)) > 0)$$
$$\leq \mathbb{N}_w(\mathcal{R}^D(W) \cap \partial D \neq \emptyset)$$
$$\leq cd(w, \partial D)^{-2}. \tag{5.4}$$

Line (5.4) follows from comparing the hitting probability with the exit probability from a ball (see Proposition 2.3 of [20]).

Combining this with (5.3) and the definition of $D_{2\epsilon}^A$ yields

$$v_\epsilon(y) \leq K\epsilon^{-2} P_y\Big(B(\tau_{2\epsilon}) \in \partial B(z, 2\epsilon)\Big).$$

We now show that for $y \in D_{4\epsilon}^A$,

$$P_y(B_{\tau_{2\epsilon}} \in \partial B(z, 2\epsilon)) \leq K P_y(B_{\tau_D} \in \Delta(z, \epsilon)). \tag{5.5}$$

This leaves us with the bound (for $d \geq 3$)

$$v_\epsilon(y) \leq K\epsilon^{-2} m_y^D(\Delta(z, \epsilon)). \tag{5.6}$$

Let $\Gamma_\alpha(\hat{n}, z)$ be the cone with opening $2\alpha$, vertex at $z$ and axis of symmetry in the $\hat{n}$ direction opening in the $\hat{n}$ direction. Since $D$ is Lipschitz we can find an $\alpha$ and a unit vector $\hat{n}$ such that $\Gamma_\alpha(\hat{n}, z)$ is an exterior cone to $D$. Let $\gamma = 1/\sin\alpha$, and set $a_\epsilon = z + \epsilon\gamma\hat{n}$. Then for $\epsilon$ sufficiently small, one has

$$B(a_\epsilon, \epsilon) \subseteq \Gamma_\alpha(\hat{n}, z) \cap D^c$$

$$B(z, 2\epsilon) \subseteq B(a_\epsilon, (\gamma + 2)\epsilon)$$
$$B(a_\epsilon, (\gamma + 3)\epsilon) \cap \partial D \subset \Delta(z, (2\gamma + 3)\epsilon).$$

As in the proof of Theorem 3.1 of [2], it follows then that
$$P_y(B_{\tau_{2\epsilon}} \in \partial B(z, 2\epsilon)) \leq K P_y(B_{\tau_D} \in \Delta(z, 5\epsilon)).$$

Finally (5.5) follows by the doubling lemma on harmonic measure (lemma 5.9 of [17], which requires $d \geq 3$).

Next we show a similar estimate holds for a lower bound, again using the same proof as in [2]. By the Cauchy-Schwarz lemma we have
$$\begin{aligned} v_\epsilon(x) &= \mathbb{N}_x(X^D(\Delta(z, \epsilon)) > 0) \\ &\geq \frac{\mathbb{N}_x(X^D(\Delta(z, \epsilon)))^2}{\mathbb{N}_x(X^D(\Delta(z, \epsilon))^2)} \\ &= \frac{m_x^D(\Delta(z, \epsilon))^2}{\int G_d(x, y) m_y^D(\Delta(z, \epsilon))^2 dy}. \end{aligned}$$

The last line following from the Palm formula, Lemma 2.6.

The lower estimate
$$v_\epsilon(x) \geq K \epsilon^{-2} m_x^D(\Delta(z, \epsilon)) \tag{5.7}$$
will follow if we show that
$$\int G_D(x, y) m_y^D(\Delta(z, \epsilon))^2 dy \leq K \epsilon^2 m_x^D(\Delta(z, \epsilon)). \tag{5.8}$$

To do this in the $d \geq 5$ case we use two results for Lipschitz domains. First, the "3-g" theorem of Cranston, Fabes and Zhao (cf. [5]) shows that there exists a constant $c$, depending on the domain $D$, for which for all $x, y, z \in D$ the following bound holds
$$\frac{G_D(x, y) G_D(y, z)}{G_D(x, z)} \leq c \left( |x - y|^{2-d} + |y - z|^{2-d} \right).$$

Second, we use the following comparison between the harmonic measure and the Green function (cf. [17], Lemma 5.8)
$$M^{-1} \leq \frac{m_y^D(\Delta(z, \epsilon))}{\epsilon^{d-2} G(y, A_\epsilon)} \leq M, \tag{5.9}$$
where $M \geq 1$ is a constant depending on $D$, $y \in D_{2\epsilon}^A$ and $A_\epsilon$ is a point satisfying $\epsilon M^{-1} \leq |A_\epsilon - z| \leq \epsilon M$ and $\epsilon M^{-1} \leq d(A_\epsilon, \partial D)$ ($A_\epsilon$ may be taken to be $z + \epsilon \hat{n}$ where $\Gamma_\alpha(\hat{n}, z)$ is an interior cone at $z$).

First, set $B = B(z, 2M\epsilon)$. Then on $B \cap D$ one has for $\epsilon$ sufficiently small that
$$G(x, y)/G(x, A_\epsilon) < K.$$

Since $m_x^D(\Delta(z, \epsilon))$ must be $\leq 1$, it follows that
$$\int_{B \cap D} G_D(x, y) m_y^D(\Delta(z, \epsilon))^2 dy \leq K G(x, A_\epsilon) \int_B dy \leq K G(x, A_\epsilon) \epsilon^d \leq K \epsilon^2 m_x^D(\Delta(z, \epsilon)).$$

The last inequality is a consequence of (5.9).

Next,

$$\int_{B^c \cap D} G_D(x,y) m_y^D(\Delta(z,\epsilon))^2 dy$$

$$\leq K \int_{B^c \cap D} \epsilon^{d-2} m_x^D(\Delta(z,\epsilon)) \frac{G_D(x,y) G_D(y, A_\epsilon)}{G_D(x, A_\epsilon)} G_D(y, A_\epsilon) dy \tag{5.10}$$

$$\leq K \epsilon^{d-2} m_x^D(\Delta(z,\epsilon)) \int_{B^c \cap D} \left( |x-y|^{2-d} + |y - A_\epsilon|^{2-d} \right) G_D(y, A_\epsilon) dy \tag{5.11}$$

$$\leq K \epsilon^{d-2} m_x^D(\Delta(z,\epsilon)) \left( K + \int_{B(A_\epsilon, M\epsilon)^c} |y - A_\epsilon|^{2(2-d)} dy \right)$$

$$\leq K \epsilon^{d-2} m_x^D(\Delta(z,\epsilon)) \left( K + \int_{M\epsilon}^\infty dr \, r^{3-d} \right)$$

$$\leq K \epsilon^2 m_x^D(\Delta(z,\epsilon)). \tag{5.12}$$

Here, (5.10) holds by (5.9), and (5.11) by the "3-g" theorem. Note that (5.12) requires that $d \geq 5$. Therefore in that case we have by (5.7), (5.6) and the boundary Harnack principle, that

$$\frac{v_\epsilon(y)}{v_\epsilon(x)} \leq K \frac{\epsilon^2 m_y^D(\Delta(z,\epsilon))}{\epsilon^2 m_x^D(\Delta(z,\epsilon))} \leq K K_x(y, z).$$

In the case $d = 4$, the argument proceeds to (5.10) just as before. Choose $\eta > 0$ and a closed circular cone $C$ with vertex $z$, such that $C \cap B(z, \eta) \subset D^c$. Set $D_0 = \mathbb{R}^4 \setminus C$. Then

$$\int_{B^c \cap D} G_D(x,y) m_y^D(\Delta(z,\epsilon))^2 dy \leq K \epsilon^2 m_x^D(\Delta(z,\epsilon)) \left( K + \int_{B^c \cap D_0} |y - A_\epsilon|^{-2} G_{D_0}(y, A_\epsilon) dy \right).$$

Thus, to show (5.8) in the case $d = 4$, it will suffice to show that

$$\int_{B^c \cap D_0} |y - A_\epsilon|^{-2} G_{D_0}(y, A_\epsilon) dy \leq K. \tag{5.13}$$

We may assume that $z = 0$. If $y, y' \in D_0$ and $|y'| \leq 1$, $|y| \geq M$, then $G_{D_0}(y, y')$ is dominated by a harmonic function of the form $|y|^{-\alpha} \Theta(|y|^{-1} y)$. Separating variables in Laplace's equation gives that $\alpha = 1 + \sqrt{1 + \beta}$, for some $\beta > 0$ depending on the opening angle of $C$. Thus $\alpha > 2$, and so by the Brownian scaling and a change of variable,

$$\int_{B^c \cap D_0} |y - A_\epsilon|^{-2} G_{D_0}(y, A_\epsilon) dy = (\epsilon^{-2})^2 \int_{B(0,M)^c \cap D_0} |y - \epsilon^{-1} A_\epsilon|^{-2} G_{D_0}(y, \epsilon^{-1} A_\epsilon) \epsilon^4 dy$$

$$\leq K \int_{B(0,M)^c \cap D_0} |y|^{-2} G_{D_0}(y, \epsilon^{-1} A_\epsilon) dy$$

$$\leq K \int_M^\infty r^{-2-\alpha} r^3 \, dr < \infty.$$

This shows (5.13).

Suppose now the lemma holds for all proper subsets of $A$. Set $\alpha = 4|A|\epsilon$. From Theorem 5.2 we have
$$\frac{1}{2}\Delta v_\epsilon^A = k_\epsilon^A v_\epsilon^A - 2 \sum_{\substack{B,C \subseteq A \\ \emptyset, A \neq B, C \\ B \cup C = A}} (-1)^{|B \cap C|} v_\epsilon^B v_\epsilon^C,$$
with $k_\epsilon^A = 2(2u_\epsilon^A + (-1)^{|A|} v_\epsilon^A) > 0$. Set
$$E = \{(B,C) : \emptyset, A \neq B, C; B \cup C = A; |B \cap C| \text{ is even}\},$$
$$O = \{(B,C) : \emptyset, A \neq B, C; B \cup C = A; |B \cap C| \text{ is odd}\}.$$

By the Feynman-Kac formula we have
$$v_\epsilon^A(y) = E_y\left(v_\epsilon^A(B_{\tau_\alpha}) \exp -\int_0^{\tau_\alpha} k_\epsilon^A(B_r)dr\right)$$
$$+ 2\sum_E E_y\left(\int_0^{\tau_\alpha}(v_\epsilon^B v_\epsilon^C)(B_t)\exp(-\int_0^t k_\epsilon^A(B_r)dr)dt\right)$$
$$- 2\sum_O E_y\left(\int_0^{\tau_\alpha}(v_\epsilon^B v_\epsilon^C)(B_t)\exp(-\int_0^t k_\epsilon^A(B_r)dr)dt\right)$$
$$\leq E_y\left(v_\epsilon^A(B_{\tau_\alpha})\exp(-\int_0^{\tau_\alpha} k_\epsilon^A(B_r)dr)\right)$$
$$+ 2\sum_E E_y\left(\int_0^{\tau_\alpha}(v_\epsilon^B v_\epsilon^C)(B_t)\exp(-\int_0^t k_\epsilon^A(B_r)dr)dt\right)$$
$$= I + II.$$

First we analyze the harmonic term $I$. Since $k_\epsilon^A \geq 0$ we have by (5.6) that
$$I \leq \sum_{i \in A} E_y(v_\epsilon^A(B_{\tau_\alpha}) \mathbf{1}(B_{\tau_\alpha} \in \partial B(z_i, \alpha)))$$
$$\leq K \sum_{i \in A} \sup_{\partial B(z_i,\alpha)} v_\epsilon^{A \setminus \{i\}}(\cdot) m_y^D(\Delta(z_i, \alpha)).$$

In particular we have
$$\frac{I}{\prod_{j \in A} v_\epsilon^j(x)} \leq K \sum_{i \in A} \left(\frac{\sup_{\partial B(z_i,\alpha)} v_\epsilon^{A \setminus \{i\}}(\cdot)}{\prod_{j \in A \setminus \{i\}} v_\epsilon^j(x)}\right) \frac{m_y^D(\Delta(z_i, \alpha))}{\epsilon^{-2} m_x^D(\Delta(z_i, \alpha))}$$
$$\leq K \sum_{i \in A} \epsilon^2 K_x(y, z_i), \tag{5.14}$$

since the first factors are bounded by induction. We see that the contribution from the harmonic term vanishes in the limit as $\epsilon \to 0$ when $|A| \geq 2$.

To bound $II$ we first remark that $v_\epsilon^B v_\epsilon^C \leq v_\epsilon^{B \setminus C} v_\epsilon^C$ (or $v_\epsilon^B v_\epsilon^{C \setminus B}$ if $B \setminus C$ is empty), and so it suffices to assume that we can bound the terms of $II$ when $B \cap C = \emptyset$. The induction hypothesis applies to proper subsets of $B \subset A$ as $\alpha > 4|A|\epsilon \geq 4|B|\epsilon$. For $B \cap C = \emptyset$, $B \cup C = A$, $B \neq \emptyset$, $C \neq \emptyset$, we have that

$$E_y\left(\int_0^{\tau_\alpha}(v_\epsilon^B v_\epsilon^C)(B_t)\exp(-\int_0^t k_\epsilon^A(B_r)dr)dt\right)$$

$$\leq K(\prod_{i\in A}v_\epsilon^i(x))E_y\left(\int_0^{\tau_D}dt\,(\sum_B K_x^D(B_t,z_i))(\sum_C K_x^D(B_t,z_j))\right).$$

As $II$ is bounded by a finite constant times terms as above, it suffices to consider terms of the following type with $i\neq j$

$$E_y\left(\int_0^{\tau_D}dt\,K_x^D(B_t,z_i)K_x^D(B_t,z_j)\right)=G_D\left(K_x^D(\cdot,z_i)K_x^D(\cdot,z_j)\right)(y). \tag{5.15}$$

Since $z_i$ and $z_j$ are separated by at least $\delta$ we have that there exists a constant $K$ depending only on $\delta$ and $D$ such that

$$G_D\left(K_x^D(\cdot,z_i)K_x^D(\cdot,z_j)\right)(y)\leq K\left(G_D(K_x^D(\cdot,z_i))(y)+G_D(K_x^D(\cdot,z_j))(y)\right).$$

The "3-g" theorem gives

$$G_D(K_x^D(\cdot,z_i))(y)=\int_D \lim_{z\to z_i}\frac{G_D(y,w)G_D(w,z)}{G_D(x,z)}dw$$

$$\leq K\int_D \lim_{z\to z_i}\frac{G_D(y,z)}{G_D(x,z)}\left(|y-w|^{2-d}+|w-z|^{2-d}\right)dw$$

$$\leq KK_x^D(y,z_i). \tag{5.16}$$

Putting together lines (5.16), (5.15), and (5.14) we get the necessary bounds for the proof of the lemma. □

Next, to finish the proof of Theorem 5.3, we again use induction on the size of $A$. First consider the case where we hit a single point. Let $A=\{i\}$. Let $\delta>0$, and recall that $\tau_\delta=\tau_\delta^A=\tau_{D_\delta^A}$. Where convenient, we will also write $D_\delta$ for $D_\delta^A$. Let $P_{yz}$ be the law of Brownian motion conditioned to leave a domain at $z$ starting from $y$.

We start with the following lemma

**Lemma 5.5.** *Let $y\in D_\delta^A$. Then uniformly in $z\in\partial D$*

$$\lim_{\epsilon\to 0}P_{yz}(\exp(-\int_0^{\tau_\delta}k_\epsilon^A(B_t)dt))=1.$$

*Proof.* We have $v_\epsilon^i(x)\to 0$, so by Lemma 5.4 it suffices to prove that

$$\lim_{\lambda\to 0}P_{yz}(\exp(-\lambda\tau_\delta))=1, \tag{5.17}$$

uniformly in $z$. But

$$\sup_{y\in D_\delta,z\in\partial D}P_{yz}(\tau_\delta)<\infty$$

since $D$ is Lipschitz [5]. Line (5.17) then follows immediately by a well-known argument of Khasminsky. See Lemma 3.7 of [4]. □

Let $x,y\in D_{\delta_0}$, by the Feynman-Kac formula, for each fixed $\delta<\delta_0$ we have

$$\frac{v_\epsilon^i(y)}{v_\epsilon^i(x)}=\frac{E_y(v_\epsilon^i(B_{\tau_\delta})\exp-\int_0^{\tau_\delta}2k_\epsilon^A(B_r)dr)}{E_x(v_\epsilon^i(B_{\tau_\delta})\exp-\int_0^{\tau_\delta}2k_\epsilon^A(B_r)dr)}$$

$$= \left( \frac{\int_{\partial D_\delta} P_{yz}(\exp - \int_0^{\tau_\delta} k_\epsilon^A(B_r)dr) K_x^{D_\delta}(y,z) v_\epsilon^i(z) m_x^{D_\delta}(dz)}{\int_{\partial D_\delta} v_\epsilon^i(z) m_x^{D_\delta}(dz)} \right)$$
$$\times \left( \frac{\int_{\partial D_\delta} v_\epsilon^i(z) m_x^{D_\delta}(dz)}{\int_{\partial D_\delta} P_{xz}(\exp - \int_0^{\tau_\delta} k_\epsilon^A(B_r)dr) v_\epsilon^i(z) m_x^{D_\delta}(dz)} \right).$$

The measure
$$\lambda_{\epsilon,\delta}(x,dz) = \frac{v_\epsilon^i(z) m_x^{D_\delta}(dz)}{\int_{\partial D_\delta} v_\epsilon^i(z) m_x^{D_\delta}(dz)} \in \mathcal{M}_1(\partial D_\delta^A).$$

Since the boundary of $D_\delta^A$ is compact, by Prohorov's theorem any sequence $\epsilon_j$ has a subsequence, again written $\epsilon_j$, for which $\lambda_{\epsilon_j,\delta} \Rightarrow \lambda_\delta \in \mathcal{M}_1(\partial D_\delta^A)$. Also, $K_x^{D_\delta}(y,z)$ is continuous and bounded in $z$, for $z \in D \cap \partial D_\delta^A$, when $x, y \in D_{\delta_0}^A$. Since $P_{yz}(\exp - \int_0^{\tau_\delta} k_\epsilon^A(B_r)dr) \to 1$ uniformly as $\epsilon \to 0$ by Lemma 5.5, we have that for $x, y \in D_{\delta_0}^A$ and for all $\delta < \delta_0$
$$\lim_{j \to \infty} \frac{v_\epsilon^i(y)}{v_\epsilon^i(x)} = \Psi_x^{(\epsilon_j)}(y) = \int_{\partial D_\delta} K_x^{D_\delta}(y,z) \lambda_\delta(x,dz).$$

This limit is then harmonic in $y$ for $y \in D_{k_0}$. By a diagonalization argument, we can assume there exists a convergent subsequence of our sequence such that the convergence holds simultaneously for a sequence of $\delta_j$'s which converge to 0. By Lemma 5.4 we see then that the limit is harmonic in $y$ with boundary value 0 on $\partial D \cap \partial D_\delta^A$ for all $\delta > 0$, and is 1 at $y = x$. In other words, $\lim_j \Psi_x^{(\epsilon_j)}(y)$ is the Martin kernel for Brownian motion in $D$. Thus all subsequences have a subsequence which converges to the Martin kernel, and so the limit itself exists.

To prove the induction step, fix $A$ and assume the result is true for all proper subsets of $A$. Therefore, if $B, C \subseteq A$ and $B, C \neq A$ we have
$$\lim_{\epsilon \to 0} \frac{v_\epsilon^B(y) v_\epsilon^C(y)}{\prod_A v_\epsilon^i(x)} = (\Psi_x^B \Psi_x^C)(y) \mathbf{1}_{B \cap C = \emptyset}.$$

Again by the Feynman-Kac formula and Theorem 5.2, we have a $\delta_0$ such that if $\delta < \delta_0$ and $x, y \in D_{\delta_0}$ then
$$\frac{v_\epsilon^A(y)}{\prod_A v_\epsilon^i(x)} \geq 2 \left( \sum_{(B,C) \in E} - \sum_{(B,C) \in O} \right) E_y \left( \int_0^{\tau_\delta} \frac{v_\epsilon^B(B_t) v_\epsilon^C(B_t)}{\prod_A v_\epsilon^i(x)} \exp(-\int_0^t k_\epsilon^A(B_r)dr) dt \right).$$

Consider the right-hand side as $\epsilon \to 0$. By Lemmas 5.4 and 5.5, the dominated convergence theorem applies, and thus by the induction hypothesis,
$$\liminf_{\epsilon \to 0} \frac{v_\epsilon^A(y)}{\prod_A v_\epsilon^i(x)} \geq 2 \sum_{\substack{B \cap C = \emptyset \\ B \cup C = A \\ B, C \neq \emptyset}} E_y\left( \int_0^{\tau_D} (\Psi_x^B \Psi_x^C)(B_t) dt \right).$$

For the upper bound, set $\alpha = 12|A|\epsilon$. Then as in the proof of Lemma 5.4, for $y \in D_\alpha^A$ the differential equation for $v_\epsilon^A$ again yields the bound
$$\frac{v_\epsilon^A(y)}{\prod_{i \in A} v_\epsilon^i(x)} \leq \frac{E_y(v_\epsilon^A(B_{\tau_\alpha}))}{\prod_{i \in A} v_\epsilon^i(x)} + 2 \sum_{(B,C) \in E} \frac{E_y(\int_0^{\tau_\alpha} (v_x^B v_x^C)(B_t) dt)}{\prod_A v_\epsilon^i(x)}.$$

As in (5.14), the first term goes to 0. Therefore we can apply the dominated convergence theorem to the second term as before, and get by induction that

$$\limsup_{\epsilon \to 0} \frac{v_\epsilon^A(y)}{\prod_A v_\epsilon^i(x)} \leq 2 \sum_{\substack{B \cap C = \emptyset \\ B \cup C = A \\ B, C \neq \emptyset}} E_y\left(\int_0^{\tau_D} (\Psi_x^B \Psi_x^C)(B_t) dt\right).$$

□

**Theorem 5.6.** *Let $D$ satisfy the conditions of Theorem 5.3, $\Phi_k \in \mathcal{F}_k$ and fix $x \in D$. Then for $y \in D$*

$$\lim_{\epsilon \to 0} \mathbb{N}_y(\Phi_k \mid \prod_{i=1}^n \langle X^D, \mathbf{1}_{\Delta_\epsilon^i} \rangle > 0) = \frac{1}{\Psi_x^N(y)} \mathbb{N}_y(\Phi_k M_k^N),$$

*where*

$$M_k^N = \sum_{\sigma \in \mathcal{P}(N)} \prod_{C \in \sigma} \langle X^k, \Psi_x^C \rangle.$$

*Remark* 5.7. By Theorem 5.3, the probabilistic representation of section 4 applies to the limiting measures

$$\mathbb{M}_y(\Phi_k) = \frac{1}{\Psi_x^N(y)} \mathbb{N}_y(\Phi_k M_k^N).$$

We take $g = 0$. The backbone is a tree connecting $y$ to the $n$ points $\{z_1, \ldots z_n\}$, which throws off mass that evolves as an unconditioned superprocess.

*Proof.* The idea is to use the special Markov property of the Exit measure to get a formal expression for the limit and then show that the convergence indeed holds. Set $W_\epsilon^C = \exp\left((-1)^{|C|}\langle X^k, v_\epsilon^C \rangle\right) - 1$, so that $W_\epsilon^C = (-1)^{|C|}|W_\epsilon^C|$. First we consider

$$\mathbb{N}_y\left(\prod_{i=1}^n \langle X^D, \mathbf{1}_{\Delta_\epsilon^i} \rangle > 0 \mid \mathcal{F}_k\right)$$

$$= \lim_{\lambda \to \infty} \mathbb{N}_y\left(\prod_{i=1}^n (1 - \exp(-\lambda \langle X^D, \mathbf{1}_{\Delta_\epsilon^i} \rangle)) \mid \mathcal{F}_k\right)$$

$$= \lim_{\lambda \to \infty} \mathbb{N}_y\left(1 + \sum_{\emptyset \neq B \subseteq N} (-1)^{|B|} \exp -\lambda \langle X^D, \sum_{i \in B} \mathbf{1}_{\Delta_\epsilon^i} \rangle \mid \mathcal{F}_k\right) \tag{5.18}$$

$$= \lim_{\lambda \to \infty} 1 + \sum_{\emptyset \neq B \subseteq N} (-1)^{|B|} \exp -\langle X^k, \mathbb{N}.(1 - \exp -\lambda \langle X^D, \sum_{i \in B} \mathbf{1}_{\Delta_\epsilon^i} \rangle) \rangle \tag{5.19}$$

$$= 1 + \sum_{\emptyset \neq B \subseteq N} (-1)^{|B|} \exp -\langle X^k, u_\epsilon^B \rangle \tag{5.20}$$

$$= 1 + \sum_{\emptyset \neq B \subseteq N} (-1)^{|B|} \exp - \sum_{\emptyset \neq C \subseteq B} (-1)^{|C|+1} \langle X^k, v_\epsilon^C \rangle \tag{5.21}$$

$$= 1 + \sum_{\emptyset \neq B \subseteq N} (-1)^{|B|} \left(1 + \sum_{j=1}^{2^{|B|}-1} \sum_{\substack{C_1 \prec C_2 \prec \cdots \prec C_j \\ \emptyset \neq C_i \subseteq B \, \forall i}} \prod_{i=1}^j W_\epsilon^{C_i}\right) \tag{5.22}$$

$$= 1 + \sum_{\emptyset \neq B \subseteq N} (-1)^{|B|} + \sum_{j=1}^{2^{|N|}-1} \sum_{\substack{C_1 \prec C_2 \prec \cdots \prec C_j \\ \emptyset \neq C_i \subseteq N \, \forall i}} \prod_{i=1}^{j} W_\epsilon^{C_i} \left( \sum_{C_1 \cup \cdots \cup C_j \subseteq B \subseteq N} (-1)^{|B|} \right)$$

$$= \sum_{j=1}^{2^{|N|}-1} \sum_{\substack{C_1 \prec C_2 \prec \cdots \prec C_j \\ C_1 \cup \cdots \cup C_j = N \\ \emptyset \neq C_i \, \forall i}} (-1)^n \prod_{i=1}^{j} W_\epsilon^{C_i} \tag{5.23}$$

$$= \sum_{j=1}^{2^{|N|}-1} \sum_{\substack{C_1 \prec C_2 \prec \cdots \prec C_j \\ C_1 \cup \cdots \cup C_j = N \\ \emptyset \neq C_i \, \forall i}} (-1)^n \left( \prod_{i=1}^{j} (-1)^{|C_i|} \right) \prod_{i=1}^{j} |W_\epsilon^{C_i}|$$

$$= \sum_{j=1}^{2^{|N|}-1} \sum_{\substack{C_1 \prec C_2 \prec \cdots \prec C_j \\ C_1 \cup \cdots \cup C_j = N \\ \emptyset \neq C_i \, \forall i}} \left( \prod_{i=1}^{j} |W_\epsilon^{C_i}| \right) (-1)^{n + \sum_{i=1}^{j} |C_i|}.$$

In the above, (5.18) and (5.22) follow from Lemma 2.2, (5.19) from Lemma 2.5, (5.20) from Lemma 2.3, (5.21) from Lemma 5.1, and (5.23) from Lemma 2.1.

Thus, we need to show that the following sorts of limits exist.

$$\lim_\epsilon \frac{\mathbb{N}_y(\Phi_k \prod_{i=1}^{j} |W_\epsilon^{C_i}|)}{\prod_{i=1}^{n} v_\epsilon^i(x)}. \tag{5.24}$$

To that end, for fixed $k$ we note that for $C \neq N$ there exists a constant $K < \infty$ for which $v_\epsilon^C(\cdot)/\prod_{i \in C} v_\epsilon^i(x) < K$ in $D_k$. Since $1 - e^{-x} \leq e^x - 1 \leq xe^x$ for $x > 0$ we have

$$\frac{|W_\epsilon^C|}{\prod_{i \in C} v_\epsilon^i(x)} \leq \frac{\exp\langle X^k, v_\epsilon^C \rangle - 1}{\prod_{i \in C} v_\epsilon^i(x)} \leq K \langle X^k, 1 \rangle \exp\left( \langle X^k, 1 \rangle K \prod_{i \in C} v_\epsilon^i(x) \right).$$

Moreover, by Theorem 5.3 (with dominated convergence verified by Lemma 5.4)

$$\lim_{\epsilon \to 0} \frac{\langle X^k, v_\epsilon^C \rangle}{\prod_{i \in C} v_\epsilon^i(x)} = \langle X^k, \Psi_x^C \rangle,$$

so also

$$\lim_{\epsilon \to 0} \frac{|W_\epsilon^C|}{\prod_{i \in C} v_\epsilon^i(x)} = \langle X^k, \Psi_x^C \rangle.$$

Set $m_\epsilon^n = \sup_{C \subsetneq N} \prod_{i \in C} v_\epsilon^i(x)$, so $m_\epsilon^n \to 0$, as $\epsilon \to 0$. This gives that for bounded $\Phi_k \in \mathcal{F}_k$,

$$\left| \frac{\Phi_k \prod_{i=1}^{j} |W_\epsilon^{C_i}|}{\prod_{i=1}^{n} v_\epsilon^i(x)} \right| \leq cK^n \langle X^k, 1 \rangle^n \exp(\langle X^k, 1 \rangle K n m_\epsilon^n) \prod_{i=1}^{n} (v_\epsilon^i(x))^{\alpha_i - 1},$$

where $\alpha_i$ is the number of times $i$ appears in the $\{C_k\}$. By Lemma 2.7, and the fact that $x^n e^x \leq K(e^{2x} - 1)$ for $x \geq 0$, this bound is in $\mathcal{L}^1(\mathbb{N}_x)$ and it decreases as $\epsilon \to 0$. Thus we may apply the dominated convergence theorem to calculate limits like (5.24). The limit is 0 unless all the $\alpha_i = 1$, so if $\{C_i\}_{i=1}^j$ are distinct, nonempty and $\cup_{i=1}^j C_i = N$, then
$$\lim_\epsilon \frac{\mathbb{N}_y(\Phi_k \prod_{i=1}^j |W_\epsilon^{C_i}|)}{\prod_{i=1}^n v_\epsilon^i(x)} = \mathbb{N}_y(\Phi_k \prod_{i=1}^j \langle X^k, \Psi_x^{C_i}\rangle) \mathbf{1}_{\{C_i \text{ disjoint}\}}.$$
The formula for $M_k^N$ follows immediately, using Lemma 3.4. □

*Remark* 5.8. Recall that Theorem 3.1 includes the hypothesis that $v^N < \infty$. With $z_1, \ldots, z_n$ distinct and $D$ Lipschitz, this condition follows from Lemma 5.4. It is worth observing that it may fail if the $z_i$ are not distinct. For example, if $D = B(0, 1)$ and $z \in \partial D$, we take $n = 2$, $z_1 = z_2 = z$, and $v^1 = v^2 = K_0(\cdot, z)$. Then for $C$ a suitable cone, with vertex at $z$,
$$v^N(x) = 2 \int_D G_D(x, y) K_0^D(y, z)^2 \, dy$$
$$\geq K \int_{C \cap B(z, 1/2)} d(y, \partial D)[d(y, \partial D)^{-(d-1)}]^2 \, dy$$
$$\geq K \int_0^{1/2} r^{1-2(d-1)} r^{d-1} \, dr$$
$$= K \int_0^{1/2} r^{2-d} \, dr = \infty,$$
for $d \geq 3$.

*Remark* 5.9. What if, instead of using the excursion measure $\mathbb{N}_y$, we use the probability measure $\mathbb{P}_\mu$, under which the superprocess $X_t$ has $X_0 = \mu$? In this case, one expects the form of the answer to change, as it is no longer necessary that the backbone originate with a single path. In genealogical terms, there may be more than one ancestor for the particles reaching the set $\{z_1, \ldots, z_n\}$.

In fact, the argument of Theorem 5.6 still applies (assuming for simplicity that $\mu$ has compact support in $D$), so that for $\Phi_k \in \mathcal{F}_k$,
$$\lim_{\epsilon \to 0} \mathbb{P}_\mu(\Phi_k \mid \prod_{i=1}^n \langle X^D, 1_{\Delta_\epsilon^i}\rangle > 0) = \frac{1}{\mathbb{P}_\mu(M_k^N)} \mathbb{P}_\mu(\Phi_k M_k^N), \tag{5.25}$$
where $M_k^N = \sum_{\sigma \in \mathcal{P}(N)} \prod_{C \in \sigma} \langle X^k, \Psi_x^C\rangle$ as before. For simplicity, drop the subscript $x$. It remains to evaluate the expressions in (5.25).

Recall that
$$\mathbb{P}_\mu(\exp -\langle X^k, \phi\rangle) = \exp -\langle \mu, \mathbb{N}.(1 - e_\phi^k)\rangle.$$
Let $\sigma \in \mathcal{P}(N)$, and take $s = |\sigma|$. As in Lemma 2.6,
$$\mathbb{P}_\mu(e_\phi^k \prod_{i=1}^s \langle X^k, \Psi^{\sigma(i)}\rangle) = (-1)^s \frac{d^s}{d\lambda_1 \ldots d\lambda_s} \mathbb{P}_\mu(e^k(\phi + \sum_{i=1}^s \lambda_i \Psi^{\sigma(i)}))$$
$$= (-1)^s \frac{d^s}{d\lambda_1 \ldots d\lambda_s} \exp -\langle \mu, \mathbb{N}.(1 - e^k(\phi + \sum_{i=1}^s \lambda_i \Psi^{\sigma(i)}))\rangle$$
$$= \mathbb{P}_\mu(e_\phi^k) \sum_{\gamma \in \mathcal{P}(\{1,\ldots,s\})} \prod_{i=1}^{|\gamma|} \langle \mu, \mathbb{N}.(e_\phi^k \prod_{j \in \gamma(i)} \langle X^k, \Psi^{\sigma(j)}\rangle)\rangle.$$

Thus

$$\mathbb{P}_\mu(e_\phi^k M_k^N) = \mathbb{P}_\mu(e_\phi^k) \sum_{\substack{\gamma \in \mathcal{P}(N) \\ \sigma_i \in \mathcal{P}(\gamma(i)) \\ i=1,\ldots,|\gamma|}} \prod_{i=1}^{|\gamma|} \langle \mu, \mathbb{N}.(e_\phi^k \prod_{A \in \sigma_i} \langle X^k, \Psi^A \rangle)\rangle$$

$$= \mathbb{P}_\mu(e_\phi^k) \sum_{\gamma \in \mathcal{P}(N)} \prod_{i=1}^{|\gamma|} \langle \mu, \mathbb{N}.(e_\phi^k \sum_{\sigma \in \mathcal{P}(\gamma(i))} \prod_{A \in \sigma} \langle X^k, \Psi^A \rangle)\rangle$$

$$= \mathbb{P}_\mu(e_\phi^k) \sum_{\gamma \in \mathcal{P}(N)} \prod_{i=1}^{|\gamma|} \langle \mu, \mathbb{N}.(e_\phi^k M_k^{\gamma(i)})\rangle$$

$$= \mathbb{P}_\mu(e_\phi^k) \sum_{\gamma \in \mathcal{P}(N)} \prod_{B \in \gamma} \langle \Psi^B \mu, \mathbb{M}^B.(e_\phi^k)\rangle,$$

where $\mathbb{M}_y^B(e_\phi^k) = \frac{1}{\Psi^B(y)}\mathbb{N}_y(e_\phi^k M_k^B)$. In particular,

$$\mathbb{P}_\mu(M_k^N) = \sum_{\gamma \in \mathcal{P}(N)} \prod_{B \in \gamma} \langle \mu, \Psi^B \rangle.$$

Thus, the conditioned exit measure $X^k$, with Laplace functional

$$\phi \mapsto \hat{\mathbb{P}}_\mu(e_\phi^k) = \frac{1}{\mathbb{P}_\mu(M_k^N)}\mathbb{P}_\mu(e_\phi^k M_k^N)$$

can be realized as a sum of two independent components. The first is a copy of the unconditioned exit measure, which we will denote $X_0^k$. To construct the other component, we choose a random $\gamma \in \mathcal{P}(N)$ with probability proportional to $\prod_{B \in \gamma}\langle \mu, \Psi^B \rangle$. For each $B \in \gamma$, we then independently choose a starting point $x$ with law $\frac{1}{\langle \mu, \Psi^B \rangle}\Psi^B \mu$, and then independently generate exit measures $X_B^k$ with Laplace functionals

$$\phi \mapsto \mathbb{M}_x^B(e_\phi^k).$$

We then have that

$$X^k = X_0^k + \sum_{B \in \gamma} X_B^k.$$

In other words, the backbone now consists of a branching forest, each tree of which has the same description as before, except that it targets a given subset of $N$.

*Remark* 5.10. One way to generalize 3.3 would be to seek martingales of the form

$$M_k = \sum_{j=1}^n \langle X^{k,\otimes j}, u_j \rangle, \tag{5.26}$$

where $X^{k,\otimes j} = X^k \otimes \cdots \otimes X^k$ (with $j$ factors), and $u_j(x_1, \ldots, x_j) \geq 0$ is a symmetric function of $j$ variables. Taking $u_{n+1} = 0$, the required conditions are that

$$\frac{1}{2}\Delta_{(x_i)}u_j(x_1, \ldots, x_j) = -2(j+1)u_{j+1}(x_i, x_1, x_2, \ldots, x_j),$$

for each $j$. An example would be the symmetrization of the function

$$\sum_{\sigma \in \mathcal{P}(N), |\sigma|=j} \prod_{i=1}^{j} v^{\sigma(i)}(x_i).$$

Note that the leading term $u_n$ is $n$-harmonic, in the sense that it is harmonic in each variable separately. Rather than trying to derive an analogue of Theorem 4.2 in this context, we simply note that there is an integral representation of the $n$-harmonic functions, in terms of products of the Martin kernel, of the form $h(x_1, \ldots, x_j) = \prod_{i=1}^{j} K^D(x_j, z_j)$ (see [3]). Thus a transform by $M_k$ of the above form will typically be a superposition of transforms of the type arising in Theorem 5.6. See [6] or [25] for other connections between $n$-harmonic functions and conditionings.

## 6. Conditioning on the Brownian snake exiting at exactly $n$ points

If we designate $n$ small balls on the boundary of $D$ and this time condition the exit measure to charge each one, but not to charge the complement of their union, then the following formal calculations suggest a representation of the exit measures as before. The new process will again correspond to a branching process with immigration, however the immigrated mass is conditioned to die before reaching the boundary of $D$. The importance of this type of conditioning is that these conditionings should be minimal, in some sense, and so should correspond to points of some generalized Martin boundary. Recall that $\Delta_\epsilon^A = \cup_{i \in A} \Delta_\epsilon^i$. The analogue of Theorem 5.2 is:

**Theorem 6.1.** *For $A \subset \{1, 2 \ldots, n\}$. Define*

$$u(x) = \mathbb{N}_x(\text{hit } \partial D) = \mathbb{N}_x(\langle X^D, \mathbf{1} \rangle > 0),$$
$$u_\epsilon^A(x) = \mathbb{N}_x(\text{hit } \Delta_\epsilon^A, \text{ miss } \partial D \setminus \Delta_\epsilon^A)$$
$$= \mathbb{N}_x(\langle X^D, \mathbf{1}_{\Delta_\epsilon^A} \rangle > 0, \langle X^D, \mathbf{1}_{\partial D \setminus \Delta_\epsilon^A} \rangle = 0),$$
$$v_\epsilon^A(x) = \mathbb{N}_x(\text{hit each } \Delta_\epsilon^i, \text{ miss } \partial D \setminus \Delta_\epsilon^A)$$
$$= \mathbb{N}_x(\langle X^D, \mathbf{1}_{\Delta_\epsilon^i} \rangle > 0 \ \forall i \in A, \langle X^D, \mathbf{1}_{\partial D \setminus \Delta_\epsilon^A} \rangle = 0).$$

*Then $v_\epsilon^A$ is twice differentiable in $D$ and satisfies*

$$\frac{1}{2} \Delta v_\epsilon^A = 4 u v_\epsilon^A - 2 \sum_{\substack{B \cup C = A \\ \emptyset \neq B, C}} v_\epsilon^B v_\epsilon^C = 2 \Big( 2(u - u_\epsilon^A) + v_\epsilon^A \Big) v_\epsilon^A - 2 \sum_{\substack{B \cup C = A \\ \emptyset, A \neq B, C}} v_\epsilon^B v_\epsilon^C,$$

*with boundary condition*

$$\lim_{y \to z \in \partial D} v_\epsilon^A(y) = \begin{cases} 0 & |A| \geq 2, \\ \begin{cases} 0 & z \in \partial D \setminus \overline{\Delta_\epsilon^i} \\ \infty & z \in \Delta_\epsilon^i \end{cases} & A = \{i\}. \end{cases}$$

**Conjecture .** We conjecture, in analogy to Theorem 5.3, that the limit

$$\lim_{\epsilon \to 0} \frac{v_\epsilon^A(y)}{\prod_{i \in A} v_\epsilon^i(x)} = \Gamma_x^A(y) \tag{6.1}$$

exists, where
$$\frac{1}{2}\Delta\Gamma_x^A = 4u\Gamma_x^A - 2\sum_{\substack{B\subseteq A \\ \emptyset, A\neq B}} \Gamma_x^B \Gamma_x^{A\setminus B}, \qquad (6.2)$$

and that moreover, as in (3.2),
$$\Gamma_x^A = 2\sum_{\substack{B\subseteq A \\ \emptyset, A\neq B}} U^{4u}(\Gamma_x^B \Gamma_x^{A\setminus B}), \text{ for } |A|\geq 1. \qquad (6.3)$$

Given this, let $M_k^N$ and $\mathbb{M}_x$ be as in (3.3) and (3.4), with $g(x) = u(x) = \mathbb{N}_x(\text{hit } \partial D)$ and $v^i = \Gamma_x^{\{i\}}$. We will give a formal argument, along the lines of Theorem 5.6, that
$$\lim_{\epsilon\to 0}\mathbb{N}_y(\Phi_k \mid \text{hit each } \Delta_\epsilon^i, \text{ miss } \partial D \setminus \Delta_\epsilon^A) = \frac{1}{\Gamma_x^N(y)}\mathbb{N}_y(\Phi_k M_k^n), \qquad (6.4)$$

for $\Phi_k \in \mathcal{F}_k$. Note that the particle representation of Theorem 4.2 would automatically apply.

For ease of notation, we will drop the $\epsilon$ subscripts. We introduce the notation
$$U^A = \{\langle X^D, \mathbb{1}_{\Delta_\epsilon^A}\rangle > 0, \langle X^D, \mathbb{1}_{\partial D\setminus\Delta_\epsilon^A}\rangle = 0\},$$
$$U_A = \{\langle X^D, \mathbb{1}_{\partial D\setminus\Delta_\epsilon^A}\rangle > 0\},$$
$$V^A = \{\langle X^D, \mathbb{1}_{\Delta_\epsilon^i}\rangle > 0 \,\forall i\in A, \langle X^D, \mathbb{1}_{\partial D\setminus\Delta_\epsilon^A}\rangle = 0\},$$
$$u_A(x) = \mathbb{N}_x(U_A),$$

so that $u^A(x) = \mathbb{N}_x(U^A)$ and $v^A(x) = \mathbb{N}_x(V^A)$.

**Lemma 6.2.** *We have the following relationships,*
*(a) $U^A = \bigcup_{\emptyset\neq B\subseteq A} V^B$ as a disjoint union, so $u^A(x) = \sum_{\emptyset\neq B\subseteq A} v^B(x)$.*
*(b) $\mathbb{1}_{U^A} = \mathbb{1}_U - \mathbb{1}_{U_A}$, where $U = \{\text{hit } \partial D\}$.*
*(c) $\mathbb{1}_{V^A} = \sum_{\emptyset\neq B\subseteq A}(-1)^{|A|+|B|}\mathbb{1}_{U^B}$.*

*Proof.* To prove the third statement we rearrange terms as follows
$$\sum_{\emptyset\neq B\subseteq A}(-1)^{|A|+|B|}\mathbb{1}_{U^B} = \sum_{\emptyset\neq B\subseteq A}(-1)^{|A|+|B|}\sum_{\emptyset\neq C\subseteq B}\mathbb{1}_{V^C} \quad \text{(by } (a)\text{)}$$
$$= \sum_{\emptyset\neq C\subseteq A}\mathbb{1}_{V^C}(-1)^{|A|}\sum_{C\subseteq B\subseteq A}(-1)^{|B|}$$
$$= \sum_{\emptyset\neq C\subseteq A}\mathbb{1}_{V^C}(-1)^{|A|}(-1)^{|A|}\mathbb{1}_{C=A} = \mathbb{1}_{V^A}.$$
$\square$

*Proof.* (of Theorem 6.1)

By $(a)$ and $(b)$ of Lemma 6.2 we have
$$v^A = \sum_{\emptyset\neq B\subseteq A}(-1)^{|A|+|B|}u^B = \sum_{\emptyset\neq B\subseteq A}(-1)^{|A|+|B|}(u - u_B).$$

Thus
$$\Delta v^A = \sum_{\emptyset\neq B\subseteq A}(-1)^{|A|+|B|}(\Delta u - \Delta u_B)$$

$$
\begin{aligned}
&= \sum_{\emptyset \neq B \subseteq A} (-1)^{|A|+|B|} 4(u^2 - u_B^2) \\
&= 4 \sum_{\emptyset \neq B \subseteq A} (-1)^{|A|+|B|} \left(u^2 - (u - u^B)^2\right) &(6.5)\\
&= 4 \sum_{\emptyset \neq B \subseteq A} (-1)^{|A|+|B|} \left(2u u^B - (u^B)^2\right) \\
&= 4\left(2u \sum_{\emptyset \neq B \subseteq A} (-1)^{|A|+|B|} u^B - \sum_{\emptyset \neq B \subseteq A} (-1)^{|A|+|B|} \Big(\sum_{\emptyset \neq C \subseteq B} v^C\Big)^2\right) &(6.6)\\
&= 4\left(2uv^A - \sum_{\emptyset \neq B \subseteq A} (-1)^{|A|+|B|} \sum_{\emptyset \neq C,D \subseteq B} v^C v^D\right) \\
&= 4\left(2uv^A - \sum_{\emptyset \neq C,D \subseteq A} (-1)^{|A|} v^C v^D \sum_{C \cup D \subseteq B \subseteq A} (-1)^{|B|}\right) \\
&= 4\left(2uv^A - \sum_{\substack{C \cup D = A \\ \emptyset \neq C,D}} v^C v^D\right) &(6.7)\\
&= 4\left(2uv^A - 2v^A \sum_{\emptyset \neq C \subseteq A} v^C + (v^A)^2 - \sum_{\substack{C \cup D = A \\ \emptyset, A \neq C,D}} v^C v^D\right) \\
&= 4\Big(2(u - u^A) + v^A\Big)v^A - 4 \sum_{\substack{C \cup D = A \\ \emptyset, A \neq C,D}} v^C v^D.
\end{aligned}
$$

Here, (6.5) and (6.6) follow by (b) and (a) of Lemma 6.2, and (6.7) follows by Lemma 2.1.

As for the boundary conditions, we note that by inclusion

$$\mathbb{N}_y(\text{hit each } \Delta_\epsilon^i, \text{ miss } \partial D \setminus \Delta_\epsilon^A) \leq \mathbb{N}_y(\text{hit } \Delta_\epsilon^j)$$

for each $j$. If $|A| \geq 2$, then each $z \in \partial D$ the latter tends to 0 for some $j$. The case $A = \{i\}$ follows as before. $\square$

Now, assuming the conjecture, we sketch an informal argument for (6.4).

$$
\begin{aligned}
&\mathbb{N}_x(\Phi_k, V^N) \\
&= \mathbb{N}_x(\Phi_k \sum_{\emptyset \neq B \subseteq N} (-1)^{n+|B|} \mathbf{1}_{U^B} &(6.8)\\
&= (-1)^n \mathbb{N}_x(\Phi_k \sum_{\emptyset \neq B \subseteq N} (-1)^{|B|}(\mathbf{1}_U - \mathbf{1}_{U_B})) &(6.9)\\
&= \lim_{\lambda \to \infty} (-1)^n \mathbb{N}_x(\Phi_k \sum_{\emptyset \neq B \subseteq N} (-1)^{|B|}\Big(\exp(-\lambda \langle X^D, 1\rangle) - \exp(-\lambda \langle X^D, 1_{\Delta_\epsilon^B}\rangle)\Big)) \\
&= (-1)^n \mathbb{N}_x\Big(\Phi_k \sum_{\emptyset \neq B \subseteq N} (-1)^{|B|}(\exp-\langle X^k, u_B\rangle - \exp-\langle X^k, u\rangle)\Big) &(6.10)\\
&= (-1)^n \mathbb{N}_x\Big(\Phi_k \exp(-\langle X^k, u\rangle) \sum_{\emptyset \neq B \subseteq N} (-1)^{|B|}(\exp\langle X^k, u^B\rangle - 1)\Big)
\end{aligned}
$$

- *page 31* -
**preprint:** August 28, 1998

$$= (-1)^n \mathbb{N}_x(\Phi_k \exp(-\langle X^k, u \rangle)) \Big(1 + \sum_{\emptyset \neq B \subseteq N} (-1)^{|B|} \prod_{\emptyset \neq C \subseteq B} \exp\langle X^k, v^C \rangle \Big)$$

$$= (-1)^n \mathbb{N}_x(\Phi_k \exp(-\langle X^k, u \rangle)) \Big(1 + \sum_{\emptyset \neq B \subseteq N} (-1)^{|B|} \prod_{\emptyset \neq C \subseteq B} \Big(1 + (\exp\langle X^k, v^C \rangle - 1)\Big)\Big).$$

Set $W^C = \exp\langle X^k, v^C \rangle - 1$, so by Lemma 2.2

$$\prod_{\emptyset \neq C \subseteq B}(1 + W^C) = 1 + \sum_{j \geq 1} \sum_{\substack{C_1 \prec \cdots \prec C_j \\ \emptyset \neq C_i \subseteq B \; \forall i}} \prod_{i=1}^{j} W^{C_i}.$$

Again, we note that $\sum_{\emptyset \neq B \subseteq N}(-1)^{|B|} = -1$ by Lemma 2.1. Thus we may continue the above calculation to arrive at

$$\mathbb{N}_x(\Phi_k, V^N) = (-1)^n \mathbb{N}_x(\Phi_k \exp(-\langle X^k, u \rangle)) \Big( \sum_{\emptyset \neq B \subseteq N} (-1)^{|B|} \sum_{j \geq 1} \sum_{\substack{C_1 \prec \cdots \prec C_j \\ \emptyset \neq C_i \subseteq B \; \forall i}} \prod_{i=1}^{j} W^{C_i} \Big)$$

$$= (-1)^n \mathbb{N}_x(\Phi_k \exp(-\langle X^k, u \rangle)) \Big( \sum_{j \geq 1} \sum_{\substack{C_1 \prec \cdots \prec C_j \\ \emptyset \neq C_i \subseteq N \; \forall i}} \prod_{i=1}^{j} W^{C_i} \sum_{\cup C_i \subseteq B \subseteq N} (-1)^{|B|} \Big)$$

$$= \mathbb{N}_x(\Phi_k \exp(-\langle X^k, u \rangle)) \Big( \sum_{j \geq 1} \sum_{\substack{C_1 \prec \cdots \prec C_j \\ \emptyset \neq C_i \; \forall i \\ \cup C_i = N}} \prod_{i=1}^{j} W^{C_i} \Big),$$

the last equation being a consequence of Lemma 2.1. Therefore, if we assume there are no problems with the convergence, we have that

$$\lim_{\epsilon \to 0} \mathbb{N}_x(\Phi_k \mid V^N) = \lim_{\epsilon \to 0} \frac{\mathbb{N}_x(\Phi_k, V^N)}{v^N(x)}$$

$$= \lim_{\epsilon \to 0} \mathbb{N}_x(\Phi_k \exp(-\langle X^k, u \rangle)) \Big( \sum_{j \geq 1} \sum_{\substack{C_1 \prec \cdots \prec C_j \\ \emptyset \neq C_i \; \forall i \\ \cup C_i = N}} \prod_{i=1}^{j} W^{C_i} \Big) \Big/ \frac{v_\epsilon^N(x)}{\prod_N v_\epsilon^i(x)} \prod_N v_\epsilon^i(x)$$

$$= \frac{\mathbb{N}_x(\Phi_k \exp(-\langle X^k, u \rangle)) \sum_{\sigma \in \mathcal{P}(N)} \prod_{C \in \sigma} \langle X_k, \Gamma_x^C \rangle}{\Gamma_x^N(x)}.$$

The last line follows as the terms with the $C$'s not disjoint disappear in the limit.

## 7. Historical processes

In previous sections we have, for simplicity, avoided the use of historical processes, even though this would have given us better results. In this section we remedy this problem, and sketch how those earlier results can be strengthened.

We start by adding to the material of section 2.4. Again, let $D \subset \mathbb{R}^d$ be a domain, and take bounded, smooth domains $D_k \Uparrow D$, with $X^k$ the exit measure from $D_k$, and $\tau_k(\xi)$ the exit time of $\xi$ from $D_k$. Recall that $(W_s, \zeta_s)$ denotes the Brownian snake. Let $W_s^t$ denote the stopped snake, $W_s^t(u) = W_s(t \wedge u)$. Define the historical process, of paths killed upon exiting $D_k$, to be

$$\langle H_t^k, \Phi \rangle = \int \Phi(W_s^t) \mathbf{1}_{\tau_k \leq t} dL_t(s),$$

where now $\Phi$ is a measurable function on path space.

We have repeatedly used the $\sigma$-fields $\mathcal{F}_k$ of information determined by the super Brownian paths prior to exiting $D_k$. Their formal definition is that

$$\mathcal{F}_k = \sigma\{H_t^k \mid t \geq 0\}.$$

Let

$$h_s^k(\Phi.) = \exp - \int_s^\infty \langle H_t^k, \Phi_t \rangle \, dt, \tag{7.1}$$

where $t \mapsto \Phi_t$ is positive and continuous. Since $H_\cdot^k$ is continuous under $\mathbb{N}_x$, or under any law absolutely continuous with respect to $\mathbb{N}_x$, such laws are characterized on $\mathcal{F}_k$ by the expectations of random variables $h_0^k(\Phi.)$. In particular, our approach to extending Theorem 4.2 will be to take the calculation of the $e^k(\phi) = \exp - \langle X^k, \phi \rangle$ under $\mathbb{M}_x$, and extend it to the $h_0^k(\Phi.)$.

To do so, we need to record versions of the Palm formulae in this context. A first step towards doing so is to redefine our underlying measures, allowing them to start from a path rather than a single point. If $w$ is a path, and $t \geq 0$, we let $\mathbb{N}_{w,t}$ be the excursion measure of the Brownian snake, started from the stopped path $w^t(s) = w(s \wedge t)$. In other words, under $\mathbb{N}_{w,t}$, we have that $W_0 = w^t$ and that the law of $\zeta_\cdot$ is that of a Brownian excursion above level $t$. Thus $L_u = 0$ for $u \leq t$, so that $H_u^k = 0$ for $u \leq t$. Alternatively, we can obtain $\mathbb{N}_{w,t}$ directly as the image of $\mathbb{N}_{w(t)}$ under the mapping $\xi \mapsto \hat{\xi}$, where

$$\hat{\xi}(s) = \begin{cases} w(s), & s \leq t \\ \xi(s-t), & s > t. \end{cases}$$

We write $P_{w,t}$ for the probability under which $B_s = w(s)$ for $s \leq t$, and for which $B_{t+\cdot}$ has the law of a Brownian motion started from $w(t)$. For $s < t$, let

$$\mathcal{N}_{s,t}(e^k(\phi), h_\cdot^k(\Phi.)) = \exp - \int_s^t du \, 4\mathbb{N}_{B,u}(1 - e^k(\phi)h_u^k(\Phi.)).$$

The basic Palm formula (2.3) then takes the form

$$\mathbb{N}_{w,t}\Big(\langle X^k, \phi \rangle e^k(\psi) h_t^k(\Psi.)\Big) = E_{w,t}\Big(\phi(B_{\tau_k}) \mathcal{N}_{t,\tau_k}(e^k(\psi), h_\cdot^k(\Psi.))\Big), \tag{7.2}$$

whenever $\tau_k(w) > t$. This may be proved as in Proposition 4.1 of [20], using the point of view of chapter 4 of [8].

Using (7.2) in place of (2.3), the following extended Palm formula may be proved, exactly as in Lemma 2.6.

**Lemma 7.1.** *Let $N = \{1, 2, \ldots n\}$, $n \geq 2$. Let $\{\psi_i\}$ be a family of measurable functions, and suppose that $\tau_k(w) > s$. Then*

$$\mathbb{N}_{w,s}(e^k(\phi) h_s^k(\Phi.) \prod_{i \in N} \langle X^k, \psi_i \rangle)$$

$$= \frac{1}{2} \sum_{\substack{M \subseteq N \\ \emptyset, N \neq M}} E_{w,s}\left(4 \int_s^{\tau_k} dt\, \mathcal{N}_{s,t}(e^k(\phi), h^k_\cdot(\Phi_\cdot))\right.$$

$$\left.\times \mathbb{N}_{B,t}\left(e^k(\phi) h^k_t(\Phi_\cdot) \prod_{i \in M} \langle X^k, \psi_i \rangle\right) \mathbb{N}_{B,t}\left(e^k(\phi) h^k_t(\Phi_\cdot) \prod_{i \in N \setminus M} \langle X^k, \psi_i \rangle\right)\right).$$

We now define a probability $Q_{w,s}$ under which $\Upsilon_t$ is a measure-valued process living on a branching tree. The tree starts as a single particle following $w$, but evolves as before following time $s$. In other words, $\Upsilon_t = \delta_{w(t)}$ for $t \leq s$, and $\Upsilon_{s+}$ has the same law under $Q_{w,s}$ as $\Upsilon_\cdot$ under $Q_{w(s)}$. Of course, we actually need a historical version of this construction, so we let $\dot\Upsilon_t$ be the atomic measure on path space, which gives unit mass to each path which is the history of an atom of $\Upsilon_t$. In the notation of (4.2),

$$\dot\Upsilon_t(d\omega) = \sum_{i=1}^{n_t} \delta_{x^i}(d\omega).$$

To define $\bar{\mathbb{N}}_{w,s}$, we again form a Poisson random measure $N^k(d\chi, d\mu)$ with intensity

$$\int_s^\infty dt \int 4\dot\Upsilon_t^k(d\omega) \tilde{\mathbb{N}}_{\omega,t}(H^k \in d\chi, X^k \in d\mu),$$

where $\chi$ belongs to the space of paths of measure-valued processes, and $\mu$ belongs to the space of measures on $\partial D_k$. We then realize the exit measure under $\bar{\mathbb{N}}_{w,s}$ as $X^k = \int \mu N^k(d\chi, d\mu)$, and the historical superprocess as $H^k = \int \chi N^k(d\chi, d\mu)$. In concrete terms,

$$\bar{\mathbb{N}}_{w,s}(e^k(\phi) h^k_s(\Phi)) = Q_{w,s}\left(\exp - \int_s^\infty dt\, 4\langle \dot\Upsilon^k_t, \tilde{\mathbb{N}}_{\cdot,t}(1 - e^k(\phi) h^k_t(\Phi))\rangle\right). \tag{7.3}$$

The generalization of Theorem 4.2 is then

**Theorem 7.2.** *In the notation of Theorem 4.2, if $v^N < \infty$ then $\mathbb{M}_x = \bar{\mathbb{N}}_x$ on $\mathcal{F}_k$.*

*Proof.* The result to be shown is that if $\tau_k(w) > s$, then

$$\mathbb{M}_{w,s}(e^k(\phi) h^k_s(\Phi)) = \bar{\mathbb{N}}_{w,s}(e^k(\phi) h^k_s(\Phi)),$$

where

$$\mathbb{M}_{w,s}(e^k(\phi) h^k_s(\Phi)) = \frac{1}{v^N(w(s))} \mathbb{N}_{w,s}(e^k(\phi) h^k_s(\Phi) M^N_k).$$

In fact, the argument proceeds exactly as before, with induction being used to prove a generalization of (4.6). The necessity of considering $\mathbb{M}_{w,s}$ and $\bar{\mathbb{N}}_{w,s}$, rather than just $\mathbb{M}_x$ and $\bar{\mathbb{N}}_x$, arises from the induction. □

Similar arguments can be used to show that other measure equations, proved before simply for expectations of random variables $e^k_\phi$, can be extended to arbitrary elements of $\mathcal{F}_k$. We leave the details to the interested reader.